\documentclass[mnsc,nonblindrev]{informs3_modified1}

\OneAndAHalfSpacedXI 

\usepackage{natbib}
 \bibpunct[, ]{(}{)}{,}{a}{}{,}%
 \def\bibfont{\small}%

\usepackage{amsfonts,mathrsfs,bm}
\usepackage{algpseudocode}
\usepackage{algorithm}
\usepackage{subcaption}
\usepackage{multirow}
\usepackage{epstopdf}
\usepackage{color}
\usepackage{accents}

\newcommand{\R}{\mathbb{R}}

\graphicspath{{figures/}}

\usepackage{multibib}
\TheoremsNumberedThrough     
\ECRepeatTheorems

\EquationsNumberedThrough    

\MANUSCRIPTNO{MS-0001-1922.65}

\begin{document}


\RUNAUTHOR{Hong, Huang, and Lam}

\RUNTITLE{Learning-based Robust Optimization}

\TITLE{Learning-based Robust Optimization: Procedures and Statistical Guarantees}

\ARTICLEAUTHORS{%
\AUTHOR{L. Jeff Hong}
\AFF{School of Management and School of Data Science, Fudan University, 670 Guoshun Road, Shanghai 200433, China, \EMAIL{hong\_liu@fudan.edu.cn}} %
\AUTHOR{Zhiyuan Huang}
\AFF{Department of Industrial and Operations Engineering, University of Michigan, 1205 Beal Ave., Ann Arbor, MI 48109, \EMAIL{zhyhuang@umich.edu}}
\AUTHOR{Henry Lam}
\AFF{Department of Industrial Engineering and Operations Research, Columbia University, 500 W.~120th St., New York, NY 10027 \EMAIL{henry.lam@columbia.edu}}
} 

\ABSTRACT{%
Robust optimization (RO) is a common approach to tractably obtain safeguarding solutions for optimization problems with uncertain constraints. In this paper, we study a statistical framework to integrate data into RO, based on learning a prediction set using (combinations of) geometric shapes that are compatible with established RO tools, and a simple data-splitting validation step that achieves finite-sample nonparametric statistical guarantees on feasibility. We demonstrate how our required sample size to achieve feasibility at a given confidence level is independent of the dimensions of both the decision space and the probability space governing the stochasticity, and discuss some approaches to improve the objective performances while maintaining these dimension-free statistical feasibility guarantees.

}%



\maketitle

%


\section{Introduction}
Many optimization problems in industrial applications contain uncertain parameters in constraints where the enforcement of feasibility is of importance. This paper aims to build procedures to find good-quality solutions for these problems that are tractable and statistically accurate for high-dimensional or limited data situations. 

To locate our scope of study, we consider situations where the uncertainty in the constraints is ``stochastic", and a risk-averse modeler wants the solution to be feasible ``most of the time" while not making the decision space overly conservative. One common framework to define feasibility in this context is via a chance-constrained program (CCP)
\begin{equation}
\text{minimize\ }f(x)\text{\ \ subject to\ \ }P(g(x;\xi)\in\mathcal A)\geq1-\epsilon\label{CCP}
\end{equation}
where $f(x)\in\mathbb R$ is the objective function, $x\in\mathbb R^d$ is the decision vector, $\xi\in\mathbb R^m$ is a random vector (i.e. the uncertainty) under a probability measure $P$, and $g(x;\xi):\mathbb R^d\times\mathbb R^m\to\Omega$ with $\mathcal A\subset\Omega$ for some space $\Omega$. Using existing terminology, we sometimes call $g(x;\xi)\in\mathcal A$ the safety condition, and $\epsilon$ the tolerance level that controls the violation probability of the safety condition. In this paper we will consider $g(x;\xi)\in\mathcal A$ as linear inequalities, which constitute the commonest class of CCPs. 


We will focus on settings where $\xi$ is observed via a finite amount of data, driven by the fact that in almost every application there is no exact knowledge about the uncertainty, and that data is increasingly ubiquitous. Our problem target is to find a solution feasible for \eqref{CCP} with a given statistical confidence (with respect to the data, in a frequentist sense) that has an objective value as small as possible.

First proposed by \cite{charnes1958cost}, \cite{charnes1959chance}, \cite{miller1965chance} and \cite{prekopa1970probabilistic}, the CCP framework \eqref{CCP} has been studied extensively in the stochastic programming literature (see \cite{prekopa2003probabilistic} for a thorough introduction), with applications spanning across reservoir system design (\cite{prekopa1978flood,prekopa1978serially}), cash matching (\cite{dentcheva2004dual}), wireless cooperative network (\cite{shi2015optimal}), inventory (\cite{lejeune2007efficient}) and production management (\cite{murr2000solution}). Though not always proper (notably when the uncertainty is deterministic or bounded; see e.g., \cite{ben2009robust} P.28--29), in many situations it is natural to view uncertainty as ``stochastic", and \eqref{CCP} provides a rigorous definition of feasibility under these situations. Moreover, \eqref{CCP} sets a framework to assimilate data in a way that avoids over-conservativeness by focusing on the ``majority" of the data, as we will exploit in this paper.

Our main contribution is a framework to integrate data into robust optimization (RO) as a tool to obtain high-quality solutions feasible in the sense defined by \eqref{CCP}. Instead of directly solving \eqref{CCP}, which is known to be challenging in general, RO operates by representing the uncertainty via a (deterministic) set, often known as the uncertainty set or the ambiguity set,  and enforces the safety condition to hold for any $\xi$ within it. By suitably choosing the uncertainty set, RO is well-known to be a tractable approximation to \eqref{CCP}. We will revisit these ideas by studying a procedural framework to construct an uncertainty set as a \emph{prediction set} for the data. This consists of approximating a high probability region via combinations of tractable geometric shapes compatible with RO. As a key development, we propose a simple  data-splitting scheme to determine the size of this region that ensures rigorous statistical performance. This framework is nonparametric and applies under minimal distributional requirements.  


In terms of basic statistical property, our approach satisfies a finite-sample confidence guarantee on the feasibility of the solution in which the minimum required sample size in achieving a given confidence is provably \emph{independent} of the dimensions of both the decision space and the underlying probability space. While finite-sample guarantees are also found in existing sampling-based methods, the dimension-free property of our approach makes it a suitable resort for certain high-dimensional and limited-data situations where previous methods break down. 

The above property, which may appear very strong, needs nonetheless be complemented with good approaches to curb over-conservativeness and maintain tractability. In particular, to reduce conservativeness, a prediction set should accurately trace the shape of data. On the other hand, to retain tractability, the set should be expressible in terms of basic geometric shapes compatible with RO techniques. We will present some techniques to construct uncertainty sets that balance these two aspects, while simultaneously achieve the basic statistical property. Nonetheless, we caution that theses techniques tie conservativeness to the set volume, while often times the former is more intricate and depends on the optimization setting at hand (see, e.g., \cite{lagoa2002distributionally}). Along this line, we also discuss a method to iterate the construction of uncertainty sets that incorporate updated optimality beliefs to improve the objective performance. 



 
Our approach is related to several existing methods for approximating \eqref{CCP}. Scenario generation (SG), pioneered by \cite{calafiore2005uncertain,calafiore2006scenario,campi2008exact,campi2011sampling} and independently suggested in the context of Markov decision processes by \cite{de2004constraint}, replaces the chance constraint in \eqref{CCP} with a collection of sampled constraints. Related work include also the sample average approximation (SAA) studied in \cite{luedtke2008sample,luedtke2010integer,luedtke2014branch}, which restricts the proportion of violated constraints and resembles the discarding approach in \cite{campi2011sampling}. SG provides explicit statistical guarantees on the feasibility of the obtained solution in terms of the confidence level, the tolerance level and the sample size. It directly approximates the chance-constrained optimization without the need of a set-based representation of the uncertainty, and hence allows a high geometric flexibility in the resulting set of violation and leads to less conservative solutions. However, in general, the sample size needed to achieve a given confidence grows linearly with the dimension of the decision space, which can be demanding for large-scale problems (as pointed out by, e.g., \cite{nemirovski2006convex}, P.971). Recent work reduce dependence on the decision dimension (and its interplay with the tolerance parameter) by, for instance, regularization (\cite{campi2013random}), tighter complexity results in terms of the support rank (\cite{schildbach2013randomized}), solution-dependent number of support constraints (\cite{campi2018wait}), one-off calibration schemes (\cite{care2014fast}), sequential validation (\cite{calafiore2011research,chamanbaz2016sequential,calafiore2017repetitive}), and hybrid approaches between RO and SG that translate scenario size requirements from decision to stochasticity space dimension (\cite{margellos2014road}). Among these, our proposed step to tune the set size is closest to the calibration approaches. However, instead of calibrating a solution obtained from a randomized program, we calibrate the coverage of an uncertainty set, and control conservativeness and tractability of the resulting RO through proper learning of its shape. 




A classical approach to approximating \eqref{CCP} uses safe convex approximation (SCA), by replacing the intractable chance constraint with an inner approximating convex constraint (such that a solution feasible for the latter would also be feasible for the former) (e.g., \cite{ben2000robust,nemirovski2003tractable,nemirovski2006convex}). This approach is intimately related to RO, as the approximating constraints are often equivalent to the robust counterparts (RC) of RO problems with properly chosen uncertainty sets (e.g., \cite{ben2009robust}, Chapters 2 and 4). The statistical guarantees provided by these approximations come from probabilistic deviation bounds, which often rely on the stochastic assumptions and the constraint structure on a worst-case basis (e.g., \cite{nemirovski2006convex}, \cite{ben2009robust} Chapter 10, \cite{ben1998robust,ben1999robust,el1998robust,bertsimas2004price,bertsimas2006tractable,bertsimas2004robust,chen2007robust,calafiore2006distributionally}). Thus, although the approach carries several advantages (e.g., in handling extraordinarily small tolerance levels), the utilized bounds can be restrictive to use in some cases. Moreover, most of the results apply to a single chance constraint; when the safety condition involves several constraints that need to be jointly maintained (known as a joint chance constraint), one typically needs to reduce it to individual constraints via the Bonferroni correction, which can add pessimism (there are exceptions, however; e.g., \cite{chen2010cvar}). On the other hand, these classical results in SCA and RO are capable of constructing uncertainty sets with well-chosen shapes, without directly using prediction set properties. 



We mention two other lines of work in approximating \eqref{CCP} that can blend with data. Distributionally robust optimization (DRO), an approach dated back to \cite{scarf1958min} and of growing interest in recent years (e.g., \cite{delage2010distributionally,wiesemann2014distributionally,goh2010distributionally,ben2013robust,lim2006model}), considers using a worst-case probability distribution for $\xi$ within an ambiguity set that represents partial distributional information. The two major classes of sets consist of distance-based constraints (statistical distance from a nominal distribution such as the empirical distribution; e.g., \cite{ben2013robust,wang2016likelihood}) and moment-and-support-type constraints (including moments, dispersion, covariance and/or support, e.g., \cite{delage2010distributionally,wiesemann2014distributionally,goh2010distributionally,hanasusanto2015ambiguous}, and shape and unimodality, e.g., \cite{popescu2005semidefinite,hanasusanto2015distributionally,van2016generalized,li2016ambiguous,lam2017tail}). To provide statistical feasibility guarantee, these uncertainty sets need to be properly calibrated from data, either via direct estimation or using the statistical implications from Bayesian (\cite{gupta2015near}) or empirical likelihood (\cite{lam2017empirical,duchi2016statistics,blanchet2016sample,lam2016recovering}) methods.
Another line of work takes a Monte Carlo viewpoint and uses sequential convex approximation (\cite{hong2011sequential,hu2013smooth}) that stochastically iterates the solution to a Karush-Kuhn-Tucker (KKT) point, which guarantees local optimality of the convergent solution. This approach can be applied to data-driven situations by viewing the data as Monte Carlo samples.

Finally, some recent RO-based approaches aim to utilize data more directly. For example, \cite{goldfarb2003robust} calibrate uncertainty sets using linear regression under Gaussian assumptions.  \cite{bertsimas2013data} study a tight value-at-risk bound on a single constraint and calibrate uncertainty sets via imposing a confidence region on the distributions that govern the bound. \cite{tulabandhula2014robust} study supervised prediction models to approximate uncertainty sets and suggest using sampling or relaxation to reduce to tractable problems. Our approach follows the general idea in these work in constructing uncertainty sets that cover the ``truth" with high confidence. 




The rest of this paper is organized as follows. Section \ref{sec:proposed} presents our procedural framework and statistical implications. Section \ref{sec:learn tractable} discusses some approaches to construct tight and tractable prediction sets. Section \ref{sec:numerics} reports numerical results and comparisons with existing  methods. Additional proofs, numerical results and useful existing theorems are presented in the Appendix.

\section{Basic Framework and Implications}\label{sec:proposed}


This section lays out our basic procedural framework and implications. First, consider an approximation of \eqref{CCP} via the RO:
\begin{equation}
\text{minimize\ \ }f(x)\text{\ \ subject to\ \ }g(x;\xi)\in\mathcal A\ \ \forall\ \xi\in\mathcal U\label{RO}
\end{equation}
where $\mathcal U\in\Omega$ is an uncertainty set. Obviously, for any $x$ feasible for \eqref{RO}, $\xi\in\mathcal U$ implies $g(x;\xi)\in\mathcal A$. Therefore, by choosing $\mathcal U$ that covers a $1-\epsilon$ content of $\xi$ (i.e., $\mathcal U$ satisfies $P(\xi\in\mathcal U)\geq1-\epsilon$), any $x$ feasible for \eqref{RO} must satisfy $P(g(x;\xi)\in\mathcal A)\geq P(\xi\in\mathcal U)\geq1-\epsilon$, implying that $x$ is also feasible for \eqref{CCP}. In other words,

\begin{lemma}
Any feasible solution of \eqref{RO} using a $(1-\epsilon)$-content set $\mathcal U$ is feasible for \eqref{CCP}. \label{RO approx}
\end{lemma}

Note that \cite{ben2009robust}, P.33 discussion point B points out that it is not necessary for an uncertainty set to contain most values of the stochasticity to induce probabilistic guarantees. Nonetheless, Lemma \ref{RO approx} provides a platform to utilize data structure easily and formulate concrete procedures, as we will describe.

\subsection{Learning Uncertainty Sets}\label{sec:procedure}
Assume a given i.i.d.~data set $D=\{\xi_1,\ldots,\xi_n\}$, where $\xi_i\in\mathbb R^m$ are sampled under a continuous distribution $P$. In view of Lemma \ref{RO approx}, our basic strategy is to construct $\mathcal U=\mathcal U(D)$ that is a $(1-\epsilon)$-content prediction set for $P$ with a prescribed confidence level $1-\delta$. In other words,
\begin{equation}
\mathbb P_D\left(P(\xi\in\mathcal U(D))\geq1-\epsilon\right)\geq1-\delta\label{basic}
\end{equation}
where we use the notation $\mathbb P_D(\cdot)$ to denote the probability taken with respect to the data $D$. Using such a  $\mathcal U$, any feasible solution of \eqref{RO} is feasible for \eqref{CCP} with the same confidence level $1-\delta$, i.e.,
\begin{lemma}
Any feasible solution of \eqref{RO} using $\mathcal U$ that satisfies \eqref{basic} is feasible for \eqref{CCP} with confidence $1-\delta$. \label{RO approx1}
\end{lemma}

\eqref{basic} only focuses on the feasibility guarantee for \eqref{CCP}, but does not speak much about conservativeness. To alleviate the latter issue, we judiciously choose $\mathcal U$ according to two criteria:
 \begin{enumerate}
 \item We prefer $\mathcal U$ that has a smaller volume, which leads to a larger feasible region in \eqref{RO} and hence a less conservative inner approximation to \eqref{CCP}. Note that, with a fixed $\epsilon$, a small $\mathcal U$ means a $\mathcal U$ that contains a high probability region (HPR) of $\xi$.\label{aspect1}
\item We prefer $\mathcal U$ such that $P(\xi\in\mathcal U(D))$ is close to, not just larger than, $1-\epsilon$ with confidence $1-\delta$. We also want the coverage probability $\mathbb P_D(P(\xi\in\mathcal U(D))\geq1-\epsilon)$ to be close to, not just larger than, $1-\delta$.\label{aspect2}
      \end{enumerate}

Moreover, $\mathcal U$ needs to be chosen to be compatible with tractable tools in RO. Though this tractability depends on the type of safety condition at hand and is problem-specific, the general principle is to construct $\mathcal U$ as an HPR that is expressed via a basic geometric set or a combination of them. 

The above discussion motivates us to propose a two-phase strategy in constructing $\mathcal U$. We first split the data $D$ into two groups, denoted $D_1$ and $D_2$, with sizes $n_1$ and $n_2$ respectively. Say $D_1=\{\xi_1^1,\ldots,\xi_{n_1}^1\}$ and $D_2=\{\xi_1^2,\ldots,\xi_{n_2}^2\}$. These two data groups are used as follows:
\\

\noindent\underline{Phase 1: Shape learning.} We use $D_1$ to approximate the shape of an HPR. Two common choices of tractable basic geometric shapes are:
\begin{enumerate}
\item\emph{Ellipsoid: }Set the shape as $\mathcal S=\{(\xi-\mu)'\Sigma^{-1}(\xi-\mu)\leq\rho\}$ for some $\rho>0$. The parameters can be chosen by, for instance, setting $\mu$ as the sample mean of $D_1$ and $\Sigma$ as some covariance matrix, e.g., the sample covariance matrix, diagonalized covariance matrix, or identity matrix.

\item\emph{Polytope: }Set the shape as $\mathcal S=\{\xi:a_i'\xi\leq b_i,i=1,\ldots,k\}$ where $a_i\in\mathbb R^m$ and $b_i\in\mathbb R$. For example, for low-dimensional data, this can be obtained from a convex hull (or an approximated version) of $D_1$, or alternately, of the data that leaves out $\lfloor n_1\epsilon\rfloor$ of $D_1$ that are in the ``periphery", e.g., having the smallest Tukey depth (e.g., \cite{serfling2002quantile,hallin2010multivariate}). It can also take the shape of the objective function when it is linear (a case of interest when using the self-improving strategy that we will describe later).

\end{enumerate}

We can also combine any of the above two types of geometric sets, such as:
\begin{enumerate}
\item\emph{Union of basic geometric sets: }Given a collection of polytopes or ellipsoids $\mathcal S_i$, take $\mathcal S = \bigcup_i\mathcal S_i$.

\item\emph{Intersection of basic geometric sets: }Given a collection of polytopes or ellipsoids $\mathcal S_i$, take $\mathcal S = \bigcap_i\mathcal S_i$.

\end{enumerate}

The choices of ellipsoids and polytopes are motivated from the tractability in the resulting RO, but they may not describe an HPR of $\xi$ to sufficient accuracy. Unions or intersection of these basic geometric sets provide more flexibility in tracking the HPR of $\xi$. For example, in the case of multi-modal distribution, one can group the data into several clusters (\cite{hastie2009unsupervised}), then form a union of ellipsoids over the clusters as $\mathcal S$. For non-standard distributions, one can discretize the space into boxes and take the union of boxes that contain at least some data, inspired by the ``histogram" method in the literature of minimum volume set learning (\cite{scott2006learning}). The intersection of basic sets is useful in handling segments of $\xi$ where each segment appears in a separate constraint in a joint CCP.
\\

\noindent\underline{Phase 2: Size calibration.} We use $D_2$ to calibrate the size of the uncertainty set so that it satisfies \eqref{basic} and moreover $P(\xi\in\mathcal U(D))\approx1-\epsilon$ with coverage $\approx1-\delta$. The key idea is to use quantile estimation on a ``dimension-collapsing" transformation of the data. More concretely, first express our geometric shape obtained in Phase 1 in the form $\{\xi:t(\xi)\leq s\}$, where $t(\cdot):\mathbb R^m\to\mathbb R$ is a transformation map from the space of $\xi$ to $\mathbb R$, and $s\in\mathbb R$. For the two geometric shapes we have considered above,
\begin{enumerate}
\item\emph{Ellipsoid: }We set $t(\xi)=(\xi-\mu)'\Sigma^{-1}(\xi-\mu)$. Then the $\mathcal S$ described in Phase 1 is equivalent to $\{\xi:t(\xi)\leq\rho\}$.

\item\emph{Polytope: }Find a point, say $\mu$, in $\mathcal S^\circ$, the interior of $\mathcal S$ (e.g., the Chebyshev center (\cite{boyd2004convex}) of $\mathcal S$ or the sample mean of $D_1$ if it lies in $\mathcal S^\circ$). Let $t(\xi)=\max_{i=1,\ldots,k}(a_i'(\xi-\mu))/(b_i-a_i'\mu)$ which is well-defined since $\mu\in\mathcal S^\circ$. Then the $\mathcal S$ defined in Phase 1 is equivalent to $\{\xi:t(\xi)\leq1\}$.
\end{enumerate}

For the combinations of sets, we suppose each individual geometric shape $\mathcal S_i$ in Phase 1 possesses a transformation map $t_i(\cdot)$. Then,
\begin{enumerate}
\item\emph{Union of the basic geometric sets: }We set $t(\xi)=\min_it_i(\xi)$ as the transformation map for $\bigcup_i\mathcal S_i$. This is because $\bigcup_i\{\xi:t_i(\xi)\leq s\}=\{\xi:\min_it_i(\xi)\leq s\}$.

\item\emph{Intersection of the basic geometric sets: }We set $t(\xi)=\max_it_i(\xi)$ as the transformation map for $\bigcap_i\mathcal S_i$. This is because $\bigcap_i\{\xi:t_i(\xi)\leq s\}=\{\xi:\max_it_i(\xi)\leq s\}$
\end{enumerate}

We overwrite the value of $s$ in the representation $\{\xi:t(\xi)\leq s\}$ as $t(\xi_{(i^*)}^2)$, where $t(\xi_{(1)}^2)<t(\xi_{(2)}^2)<\cdots<t(\xi_{(n_2)}^2)$ are the ranked observations of $\{t(\xi_i^2)\}_{i=1,\ldots,n_2}$, and
\begin{equation}
i^*=\min\left\{r:\sum_{k=0}^{r-1}\binom{n_2}{k}(1-\epsilon)^k\epsilon^{n_2-k}\geq1-\delta,\ 1\leq r\leq n_2\right\}\label{quantile choice}
\end{equation}
This procedure is valid if such an $i^*$ can be found, or equivalently $1-(1-\epsilon)^{n_2}\geq1-\delta$.

\subsection{Basic Statistical Guarantees}\label{sce:basic guarantee}
Phase 1 focuses on Criterion \ref{aspect1} in Section \ref{sec:procedure} by learning the shape of an HPR. Phase 2 addresses our basic requirement \eqref{basic} and Criterion \ref{aspect2}. The choice of $s$ in Phase 2 can be explained by the elementary observation that, for any arbitrary i.i.d.~data set of size $n_2$ drawn from a continuous distribution, the $i^*$-th ranked observation as defined by \eqref{quantile choice} is a valid $1-\delta$ confidence upper bound for the $1-\epsilon$  quantile of the distribution:

\begin{lemma}
Let $Y_1,\ldots,Y_{n_2}$ be i.i.d. data in $\mathbb R$ drawn from a continuous distribution. Let $Y_{(1)}<Y_{(2)}<\cdots<Y_{(n_2)}$ be the order statistics. A $1-\delta$ confidence upper bound for the $(1-\epsilon)$-quantile of the underlying distribution is $Y_{(i^*)}$, where
$$i^*=\min\left\{r:\sum_{k=0}^{r-1}\binom{n_2}{k}(1-\epsilon)^k\epsilon^{n_2-k}\geq1-\delta,\ 1\leq r\leq n_2\right\}$$
If $\sum_{k=0}^{n_2-1}\binom{n_2}{k}(1-\epsilon)^k\epsilon^{n_2-k}<1-\delta$ or equivalently $1-(1-\epsilon)^{n_2}<1-\delta$, then none of the $Y_{(r)}$'s is a valid confidence upper bound.

Similarly, a $1-\delta$ confidence lower bound for the $(1-\epsilon)$-quantile of the underlying distribution is $Y_{(i_*)}$, where
$$i_*=\max\left\{r:\sum_{k=r}^{n_2}\binom{n_2}{k}(1-\epsilon)^k\epsilon^{n_2-k}\geq1-\delta,\ 1\leq r\leq n_2\right\}$$
If $\sum_{k=1}^{n_2}\binom{n_2}{k}(1-\epsilon)^k\epsilon^{n_2-k}<1-\delta$ or equivalently $1-\epsilon^{n_2}<1-\delta$, then none of the $Y_{(r)}$'s is a valid confidence lower bound.\label{quantile estimation}
\end{lemma}

\proof{Proof of Lemma \ref{quantile estimation}.}
Let $q_{1-\epsilon}$ be the $(1-\epsilon)$-quantile, and $F(\cdot)$ and $\bar F(\cdot)$ be the distribution function and tail distribution function of $Y_i$. Consider
\begin{align*}
P(Y_{(r)}\geq q_{1-\epsilon})&=P(\text{$\leq r-1$ of the data $\{Y_1,\ldots,Y_n\}$ are $<q_{1-\epsilon}$})\\
&=\sum_{k=0}^{r-1}\binom{n_2}{k}F(q_{1-\epsilon})^k{\bar F(q_{1-\epsilon})}^{n_2-k}\\
&=\sum_{k=0}^{r-1}\binom{n_2}{k}(1-\epsilon)^k\epsilon^{n_2-k}
\end{align*}
by the definition of $q_{1-\epsilon}$. Hence any $r$ such that $\sum_{k=0}^{r-1}\binom{n_2}{k}(1-\epsilon)^k\epsilon^{n_2-k}\geq1-\delta$ is a $1-\delta$ confidence upper bound for $q_{1-\epsilon}$, and we pick the smallest one. Note that if $\sum_{k=0}^{n_2-1}\binom{n_2}{k}(1-\epsilon)^k\epsilon^{n_2-k}<1-\delta$, then none of the $Y_{(r)}$ is a valid confidence upper bound.

Similarly, we have
\begin{align*}
P(Y_{(r)}\leq q_{1-\epsilon})&=P(\text{$\geq r$ of the data $\{Y_1,\ldots,Y_n\}$ are $\leq q_{1-\epsilon}$})\\
&=\sum_{k=r}^{n_2}\binom{n_2}{k}F(q_{1-\epsilon})^k{\bar F(q_{1-\epsilon})}^{n_2-k}\\
&=\sum_{k=r}^{n_2}\binom{n_2}{k}(1-\epsilon)^k\epsilon^{n_2-k}
\end{align*}
by the definition of $q_{1-\epsilon}$. Hence any $r$ such that $\sum_{k=r}^{n_2}\binom{n_2}{k}(1-\epsilon)^k\epsilon^{n_2-k}\geq1-\delta$ will be a $1-\delta$ confidence lower bound for $q_{1-\epsilon}$, and we pick the largest one. Note that if $\sum_{k=1}^{n_2}\binom{n_2}{k}(1-\epsilon)^k\epsilon^{n_2-k}<1-\delta$, then none of the $Y_{(r)}$ is a valid confidence lower bound.\Halmos

\endproof

Similar results in the above simple order statistics calculation can be found in, e.g., \cite{serfling2009approximation} Section 2.6.1. A key element of our procedure is that $t(\cdot)$ is constructed using only Phase 1 data $D_1$, which are independent of Phase 2. Lemma \ref{quantile estimation} implies that, conditional on $D_1$, $P(t(\xi)\leq t(\xi_{(i^*)}^2))\geq1-\epsilon$ with a (conditional) confidence $1-\delta$. From this, we can average over the realizations of $D_1$ to obtain a valid coverage for the resulting uncertainty set in the sense of satisfying \eqref{basic}. This is summarized formally as:
\begin{theorem}[Basic statistical guarantee]
Suppose $D$ is an i.i.d. data set drawn from a continuous distribution $P$ on $\mathbb R^m$, and we partition $D$ into two sets $D_1=\{\xi^1_i\}_{i=1,\ldots,n_1}$ and $D_2=\{\xi^2_i\}_{i=1,\ldots,n_2}$. Suppose $n_2\geq\log\delta/\log(1-\epsilon)$. Consider the set $\mathcal U=\mathcal U(D)=\{\xi:t(\xi)\leq s\}$, where $t:\mathbb R^m\to\mathbb R$ is a map constructed from $D_1$ such that $t(\xi)$, with $\xi$ distributed according to $P$, is a continuous random variable, and $s=t(\xi^2_{(i^*)})$ is calibrated from $D_2$ with $i^*$ defined in \eqref{quantile choice}. Then $\mathcal U$ satisfies \eqref{basic}. Consequently, an optimal solution obtained from \eqref{RO} using this $\mathcal U$ is feasible for \eqref{CCP} with confidence $1-\delta$. \label{main result}
\end{theorem}

\proof{Proof of Theorem \ref{main result}.}
Since $t(\cdot)$ depends only on $D_1$ but not $D_2$, we have, conditional on any realization of $D_1$,
\begin{equation}
\mathbb P_{D_2}(P(\xi\in\mathcal U(D))\geq1-\epsilon|D_1)=\mathbb P_{D_2}(P(t(\xi)\leq t(\xi^2_{(i^*)}))\geq1-\epsilon|D_1)=\mathbb P_{D_2}(q_{1-\epsilon}\leq t(\xi^2_{(i^*)})|D_1)\geq1-\delta\label{interim}
\end{equation}
where $q_{1-\epsilon}$ is the $(1-\epsilon)$-quantile of $t(\xi)$ (which depends on $D_1$). The first equality in \eqref{interim} follows from the representation of $\mathcal U=\{\xi:t(\xi)\leq t(\xi_{(i^*)}^2)\}$, the second equality uses the definition of a quantile, and the last inequality follows from Lemma \ref{quantile estimation} using the condition $1-(1-\epsilon)^{n_2}\geq1-\delta$, or equivalently $n_2\geq\log\delta/\log(1-\epsilon)$. Note that \eqref{interim} holds given any realization of $D_1$. Thus, taking expectation with respect to $D_1$ on both sides in \eqref{interim}, we have
$$\mathbb E_{D_1}[\mathbb P_{D_2}(P(\xi\in\mathcal U(D))\geq1-\epsilon|D_1)]\geq1-\delta$$
where $\mathbb E_{D_1}[\cdot]$ denotes the expectation with respect to $D_1$, which gives
$$\mathbb P_{D}(P(\xi\in\mathcal U(D))\geq1-\epsilon)\geq1-\delta$$
We therefore arrive at \eqref{basic}. Finally, Lemma \ref{RO approx1} guarantees that an optimal solution obtained from \eqref{RO} using the constructed $\mathcal U$ is feasible for \eqref{CCP} with confidence $1-\delta$.\Halmos
\endproof

Theorem \ref{main result} implies the validity of the approach in giving a feasible solution for CCP \eqref{CCP} with confidence $1-\delta$ for any finite sample size, as long as it is large enough such that $n_2\geq\log\delta/\log(1-\epsilon)$. The reasoning of the latter restriction can be seen easily in the proof, or more apparently from the following argument: In order to get an upper confidence bound for the quantile by choosing one of the ranked statistics, we need the probability of at least one observation to upper bound the quantile to be at least $1-\delta$. In other words, we need $P\left(\text{at least one\ }t(\xi_i^2)\geq (1-\epsilon)\text{-quantile}\right)\geq1-\delta$ or equivalently $1-(1-\epsilon)^{n_2}\geq1-\delta$.

We also mention the convenient fact that, conditional on $D_1$,
\begin{equation}
P(\xi\in\mathcal U)=P(t(\xi)\leq t(\xi_{(i^*)}^2))=F(t(\xi_{(i^*)}^2))\stackrel{d}{=}U_{(i^*)}\label{order stat}
\end{equation}
where $F(\cdot)$ is the distribution function of $t(\xi)$ and $U_{(i^*)}$ is the $i^*$-th ranked variable among $n_2$ uniform variables on $[0,1]$, and ``$\stackrel{d}{=}$" denotes equality in distribution. In other words, the theoretical tolerance level induced by our constructed uncertainty set, $P(\xi\in\mathcal U)$, is distributed as the $i^*$-th order statistic of uniform random variables, or equivalently $Beta(i^*,n_2-i^*+1)$, a Beta variable with parameters $i^*$ and $n_2-i^*+1$. Note that $P(Beta(i^*,n_2-i^*+1)\geq1-\epsilon)=P(Bin(n_2,1-\epsilon)\leq i^*-1)$ where $Bin(n_2,1-\epsilon)$ denotes a binomial variable with number of trials $n_2$ and success probability $1-\epsilon$. This informs an equivalent expression of \eqref{quantile choice} as
\begin{eqnarray*}
&&\min\left\{r:P(Beta(r,n_2-r+1)\geq1-\epsilon)\geq1-\delta,\ 1\leq r\leq n_2\right\}\\
&=&\min\left\{r:P(Bin(n_2,1-\epsilon)\leq r-1)\geq1-\delta,\ 1\leq r\leq n_2\right\}
\end{eqnarray*}


To address Criterion \ref{aspect2} in Section \ref{sec:procedure}, we use the following asymptotic behavior as $n_2\to\infty$:
\begin{theorem}[Asymptotic tightness of tolerance and confidence levels]
Under the same assumptions as in Theorem \ref{main result}, we have, conditional on $D_1$:
\begin{enumerate}
\item $P(\xi\in\mathcal U)\to1-\epsilon$ in probability (with respect to $D_2$) as $n_2\to\infty$.\label{thm pt2}
\item $\mathbb P_{D_2}(P(\xi\in\mathcal U)\geq1-\epsilon|D_1)\to1-\delta$ as $n_2\to\infty$.\label{thm pt1}
\end{enumerate}\label{consistency}
\end{theorem}
Theorem \ref{consistency} confirms that $\mathcal U$ is tightly chosen in the sense that the tolerance level and the confidence level are held asymptotically exact. This can be shown by using \eqref{order stat} together with an invocation of the Berry-Essen Theorem (\cite{durrett2010probability}) applied on the normal approximation to binomial distribution. Appendix \ref{sec:proof} shows the proof details, which use techniques similar to \cite{li2008multivariate} and \cite{serfling2009approximation} Section 2.6. In fact, one could further obtain that our choice of $i^*$ satisfies $\sqrt{n_2}\left(i^*/n_2-(1-\epsilon)\right)\to\sqrt{(1-\epsilon)\epsilon}\Phi^{-1}(1-\delta)$ as $n_2\to\infty$. As a result, the theoretical tolerance level $P(\xi\in\mathcal U)$ given $D_1$ concentrates at $1-\epsilon$ by being approximately $(1-\epsilon)+Z/\sqrt{n_2}$ where $Z\sim N\left(\sqrt{\epsilon(1-\epsilon)}\Phi^{-1}(1-\delta),\epsilon(1-\epsilon)\right)$. For further details, see Appendix \ref{sec:proof}. 


Note that, because of the discrete nature of our quantile estimate, the theoretical confidence level is not a monotone function of the sample size, and neither is there a guarantee on an exact confidence level at $1-\delta$ using a finite sample (see Appendix \ref{sec:phase 2 data}). On the other hand, Theorem \ref{consistency} Part \ref{thm pt1} guarantees that asymptotically our construction can achieve an exact confidence level.




The idea of using a dimension-collapsing transformation map $t(\cdot)$ resembles the notion of data depth in the literature of generalized quantile (\cite{li2008multivariate,serfling2002quantile}). In particular, the data depth of an observation is a positive number that measures the position of the observation from the ``center" of the data set. The larger the data depth, the closer the observation is to the center. For example, the half-space depth is the minimum number of observations on one side of any line passing through the chosen observation (\cite{hodges1955bivariate,tukey1975mathematics}), and the simplicial depth is the number of simplices formed by different combinations of observations surrounding an observation (\cite{liu1990notion}). Other common data depths include the ellipsoidally defined Mahalanobis depth (\cite{mahalanobis1936generalized}) and projection-based depths (\cite{donoho1992breakdown,zuo2003projection}). Instead of measuring the position of the data relative to the center as in the data depth literature, our transformation map is constructed to create uncertainty sets with good geometric and  tractability properties. 


\subsection{Dimension-free Sample Size Requirement}\label{sec:implications}
Theorem \ref{main result} and the associated discussion above states that we need at least $n_2\geq\log\delta/\log(1-\epsilon)$ observations in Phase 2 to construct an uncertainty set that guarantees a feasible solution for \eqref{CCP} with confidence $1-\delta$. From a \emph{purely} feasibility viewpoint, this lower bound on $n_2$ is the minimum total sample size we need: Regardless of what shape we generate in Phase 1, as long as we can express it in terms of the $t(\cdot)$ and have $\log\delta/\log(1-\epsilon)$ Phase 2 observations, the basic feasibility guarantee \eqref{basic} is attained. This number does not depend on the dimension of the decision space or the probability space. It does, however, depend roughly linearly on $1/\epsilon$ for small $\epsilon$, a drawback that is also common among sampling-based approaches including both SG and SAA and gives more edge to using safe convex approximation when applicable.

We should caution, however, that if we take $n_1=0$ or choose an arbitrary shape in Phase 1, the resulting solution is likely extremely conservative in terms of objective performance. To combat this issue, it is thus recommended to set aside some data for Phase 1 with the help of established methods borrowed from statistical learning (Section \ref{sec:learn tractable} and Appendices \ref{sec:comparison} and \ref{sec:single_linear_ml} discuss these).

\subsection{Enhancing Optimality Performance via Self-improving Reconstruction}\label{sec:enhancement}
We propose a mechanism, under the framework in Section \ref{sce:basic guarantee}, to improve the performance of an uncertainty set by incorporating updated optimality belief.


\subsubsection{An Elementary Explanation}\label{sec:intuition}
As indicated at the beginning of this section, the RO we construct is a conservative approximation to the CCP. A question is whether there is an ``optimal" uncertainty set, in the sense that it is a $(1-\epsilon)$-level prediction set, and at the same time gives rise to the same solution between the RO and the CCP. As a first observation, the uncertainty set $\mathcal U=\{\xi:g(x^*;\xi)\in\mathcal A\}$, where $x^*$ is an optimal solution to the CCP, satisfies both properties: By the definition of $x^*$, this set contains $(1-\epsilon)$-content of $P$. Moreover, when we use this $\mathcal U$ in \eqref{RO}, $x^*$ is trivially a feasible solution. Since this RO is an inner approximation to CCP, $x^*$ is optimal for both the RO and the CCP. The catch, of course, is that in reality we do not know what is $x^*$. Our suggestion is to replace $x^*$ with some approximate solution $\hat x$, leading to a set $\{\xi:g(\hat x,\xi)\in\mathcal A\}$. 




Alternately, the conservativeness of the RO can be reasoned from the fact that $\xi\in\mathcal U$, independent of what the obtained solution $\hat x$ is in \eqref{RO}, implies that $g(\hat x;\xi)\in\mathcal A$. Thus our target tolerance probability $P(g(\hat x;\xi)\in\mathcal A)$ satisfies $P(g(\hat x;\xi)\in\mathcal A)\geq P(\xi\in\mathcal U)$, and, in the presence of data, makes the actual confidence level (namely $\mathbb P_D(P(g(\hat x;\xi)\in\mathcal A)\geq1-\epsilon)$) potentially over-conservative. However, this inequality becomes an equality if $\mathcal U$ is exactly $\{\xi:g(\hat x;\xi)\in\mathcal A\}$. This suggests again that, on a high level, an uncertainty set that resembles the form $g(\hat x;\xi)\in\mathcal A$ is less conservative and preferable.

Using the above intuition, a proposed strategy is as follows. Consider finding a solution for \eqref{CCP}. In Phase 1, find an approximate HPR of the data (using some suggestions in Section \ref{sec:learn tractable}) with a reasonably chosen size (e.g., just enough to cover $(1-\epsilon)$ of the data points). Solve the RO problem using this HPR to obtain an initial solution $\hat x_0$. Then reshape the uncertainty set as $\{\xi:g(\hat x_0;\xi)\in\mathcal A\}$. Finally, conduct Phase 2 by tuning the size of this reshaped set, say we get $\{\xi:g(\hat x_0;\xi)\in\tilde{\mathcal A}\}$ where $\tilde{\mathcal A}$ is size-tuned. The final RO is:
\begin{equation}
\text{minimize\ }f(x)\text{\ \ subject to\ \ }g(x,\xi)\in\mathcal A\ \ \forall\ \xi:g(\hat x_0;\xi)\in\tilde{\mathcal A}\label{general}
\end{equation}
Evidently, if the tuning step can be done properly, i.e., the set $\{\xi:g(\hat x_0;\xi)\in\mathcal A\}$ can be expressed in the form $\{\xi:t(\xi)\leq s\}$ and $s$ is calibrated using the method in Section \ref{sec:procedure}, then the procedure retains the overall statistical confidence guarantees presented in Theorems \ref{main result} and \ref{consistency}. For convenience, we call the RO  \eqref{general} created from $\hat x_0$ and the discussed procedure a ``reconstructed" RO.

More explicitly, consider the safety condition $g(x;\xi)\in\mathcal A$ in the form of linear inequalities $Ax\leq b$ where $A\in\mathbb R^{l\times d}$ is stochastic and $b\in\mathbb R^l$ is constant. After we obtain an initial solution $\hat x_0$, we set the uncertainty set as $\mathcal U=\{A:A\hat x_0\leq b+sk\}$ where $k=(k_i)_{i=1,\ldots,l}\in\mathbb R^l$ is some positive vector and $s\in\mathbb R$. The value of $s$ is calibrated by letting $t(A)=\max_{i=1,\ldots,l}\{(a_i'\hat{x}_0-b_i)/k_i\}$ where $a_i'$ is the $i$-th row of $A$ and $b_i$ is the $i$-th entry of $b$, and $s$ is chosen as $t(A_{(i^*)}^2)$, the order statistic of Phase 2 data as defined in Section \ref{sec:procedure}. Using the uncertainty set $\mathcal U$, the constraint $Ax\leq b\ \forall\ A\in\mathcal U$ becomes $\max_{a_i'\hat x_0\leq b_i+sk_i}a_i'x\leq b_i,i=1,\ldots,l$ via constraint-wise projection of the uncertainty set, which can be reformulated into linear constraints by using standard RO machinery (see, e.g., Theorem \ref{thm:polytope}).

\subsubsection{Properties of Self-improving Reconstruction}\label{sec:theory}
We formalize the discussion in Section \ref{sec:intuition} by showing some properties of the optimization problem \eqref{general}. We focus on the setting of inequalities-based safety conditions
\begin{equation}
\text{minimize\ }f(x)\text{\ \ subject to\ \ }P(g(x;\xi)\leq b)\geq1-\epsilon\label{CCP iterative}
\end{equation}
where $g(x;\xi)=(g_j(x;\xi))_{j=1,\ldots,l}\in\mathbb R^l$ and $b=(b_j)_{j=1,\ldots,l}\in\mathbb R^l$. Suppose $\hat x_0$ is a given solution (not necessarily feasible). Suppose for now that there is a way to compute quantiles exactly for functions of $\xi$, and consider the reconstructed RO
\begin{equation}
\text{minimize\ }f(x)\text{\ \ subject to\ \ }g(x,\xi)\leq b\ \ \forall\ \xi:g(\hat x_0;\xi)\leq b+\rho k\label{RO iterative}
\end{equation}
where $k=(k_j)_{j=1,\ldots,l}\in\mathbb R^l$ is a positive vector, and $\rho=\rho(\hat x_0)$ is the $(1-\epsilon)$-quantile of $\max_{j=1,\ldots,l}\{(g_j(\hat x_0;\xi)-b_j)/k_j\}$. A useful observation is:
\begin{theorem}[Feasibility guarantee for reconstruction]
Given any solution $\hat x_0$, if $\rho$ is the $(1-\epsilon)$-quantile of $\max_{j=1,\ldots,l}\{(g_j(\hat x_0;\xi)-b_j)/k_j\}$, then any feasible solution of \eqref{RO iterative} is also feasible for \eqref{CCP iterative}.\label{feasibility}
\end{theorem}

\proof{Proof of Theorem \ref{feasibility}.}
Since $\{\xi:g(\hat x_0;\xi)\leq b+\rho k\}$ is by construction a $(1-\epsilon)$-content set for $\xi$ under $P$, Lemma \ref{RO approx} concludes the theorem immediately.\Halmos
\endproof

Note that Theorem \ref{feasibility} holds regardless of whether $\hat x_0$ is feasible for \eqref{CCP iterative}. That is, \eqref{RO iterative} is a way to output a feasible solution from the input of a possibly infeasible $\hat x_0$. What is more, in the case that $\hat x_0$ is feasible, \eqref{RO iterative} is guaranteed to give a solution at least as good:

\begin{theorem}[Monotonic objective improvement]
Under the same assumption as Theorem \ref{feasibility}, an optimal solution $\hat x$ of \eqref{RO iterative} is feasible for \eqref{CCP iterative}. Moreover, if $\hat x_0$ is feasible for \eqref{CCP iterative}, then $\hat x$ satisfies $f(\hat x)\leq f(\hat x_0)$.\label{monotone}
\end{theorem}

\proof{Proof of Theorem \ref{monotone}.}
Note that if $\hat x_0$ is feasible for \eqref{CCP iterative}, we must have $\rho\leq0$ (or else the chance constraint does not hold) and hence $\hat x_0$ must be feasible for \eqref{RO iterative}. By the optimality of $\hat x$ for \eqref{RO iterative} we must have $f(\hat x)\leq f(\hat x_0)$. The theorem concludes by invoking Theorem \ref{feasibility} that implies $\hat x$ is feasible for \eqref{CCP iterative}.\Halmos
\endproof

Together, Theorems \ref{feasibility} and \ref{monotone} give a mechanism to improve any input solution in terms of either feasibility or optimality for \eqref{CCP iterative}: If $\hat x_0$ is infeasible, then \eqref{RO iterative} corrects the infeasibility and gives a feasible solution; if $\hat x_0$ is feasible, then \eqref{RO iterative} gives a feasible solution that has an objective value at least as good.

Similar statements hold if the quantile $\rho$ is only calibrated under a given statistical confidence. To link our discussion to the procedure in Section \ref{sec:procedure}, suppose that a solution $\hat x_0$ is obtained from an RO formulation (or in fact, any other procedures) using only Phase 1 data. We have:

\begin{corollary}[Feasibility guarantee for reconstruction under statistical confidence]
Given any solution $\hat x_0$ obtained using Phase 1 data, suppose $\rho$ is the upper bound of the $(1-\epsilon)$-quantile of $\max_{j=1,\ldots,l}\{(g_j(\hat x_0;\xi)-b_j)/k_j\}$ with confidence level $1-\delta$ generated under Phase 2 data. Any feasible solution of \eqref{RO iterative} is also feasible for \eqref{CCP iterative} with the same confidence.\label{feasibility confidence}
\end{corollary}

\begin{corollary}[Improvement from reconstruction under statistical confidence]
Under the same assumptions as Corollary \ref{feasibility confidence}, an optimal solution $\hat x$ of \eqref{RO iterative} is feasible for \eqref{CCP iterative} with  confidence $1-\delta$. Moreover, if $\rho\leq0$, then $\hat x$ satisfies $f(\hat x)\leq f(\hat x_0)$.\label{confidence improvement}
\end{corollary}

The proofs of Corollaries \ref{feasibility confidence} and \ref{confidence improvement} are the same as those of Theorems \ref{feasibility} and \ref{monotone}, except that Lemma \ref{RO approx1} is invoked instead of Lemma \ref{RO approx}. Note that $\rho\leq0$ in Corollary \ref{confidence improvement} implies that $\hat x_0$ is feasible for \eqref{CCP iterative} with confidence $1-\delta$. However, the case $\rho>0$ in Corollary \ref{confidence improvement} does not directly translate to a conclusion that $\hat x_0$ is infeasible under confidence $1-\delta$, since $\rho$ is a confidence upper bound, instead of lower bound, for the quantile. This implies a possibility that $\hat x_0$ is feasible and close to the boundary of the feasible region. There is no guarantee of objective improvement under the reconstructed RO in this case, but there is still guarantee that the output $\hat x$ is feasible with confidence $1-\delta$.

Our numerical experiments in Section \ref{sec:numerics} show that, when applicable, such reconstructions frequently lead to notable improvements. Nonetheless, we caution that, depending on the constraint structure, the reconstruction step does not always lead to a significant or a strict improvement even if $\rho\leq0$, and in these cases some transformation of the constraint is needed. For example, in the case of single linear chance constraint in the form \eqref{CCP iterative} with $l=1$ and a bilinear $g(x;\xi)$, the reconstructed uncertainty set consists of one linear constraint. Consequently, the dualization of the RO (see Theorem \ref{thm:polytope}) consists of one dual variable, which optimally scales $\hat x_0$ by a scalar factor. When $b$ in \eqref{CCP iterative}  (with $l=1$) is also a stochastic source, no scaling adjustment is allowed because the ``decision variable" associated with $b$ (viewing $b$ as a random coefficient in the linear constraint) is constrained to be 1. Thus, the proposed reconstruction will show no strict improvement. However, this behavior could be avoided by suitably re-expressing the constraint. When $b$ is say positively distributed (or very likely so), one can divide both sides of the inequality by $b$ to obtain an equivalent inequality with right hand side fixed to be 1. This equivalent constraint is now improvable by our reconstruction (and the new stochasticity now comprises the ratios of the original variables, which can still be observed from the data). 

\section{Constructing Uncertainty Sets}
\label{sec:learn tractable}
Our proposed strategy in Section \ref{sec:proposed} requires constructing an uncertainty set that is tractable for RO, and recommends to trace the shape of an HPR as much as possible. 
Regarding tractability, linear RO with the uncertainty set shapes mentioned in Section \ref{sec:procedure} can be reformulated into standard optimization formulations. For convenience we document some of these results in Appendix \ref{sec:RO results}, along with some explanation on how to identify $t(\cdot)$ for the size calibration in our procedure. 

Since taking unions or intersections of basic sets gives more capability to trace HPR, we highlight the following two immediate observations. First is that unions of basic sets preserve the tractability of the robust counterpart associated with each union component, with a linear growth of the number of constraints against the number of components.
\begin{lemma}[Reformulating unions of sets]
The constraint
$$g(x;\xi)\in\mathcal A\ \ \forall\ \xi\in\mathcal U$$
where $\mathcal U=\bigcup_{i=1}^k\mathcal U^i$ is equivalent to the joint constraints
$$g(x;\xi)\in\mathcal A\ \ \forall\ \xi\in\mathcal U^i,\ \ i=1,\ldots,k$$
\label{lemma:union}
\end{lemma}

Second, in the special case of intersections of sets where each intersection component is on the portion of the stochasticity associated with each of multiple constraints, the projective separability property of uncertainty sets (e.g., \cite{ben2009robust}) gives the following:
\begin{lemma}[Reformulating intersections of sets]
Let $\xi\in\mathbb R^m$ be a vector that can be represented as $\xi=(\xi^i)_{i=1,\ldots,k}$, where $\xi^i\in\mathbb R^{m^i},i=1,\ldots,k$ are vectors such that $\sum_{i=1}^km^i=m$. Suppose that $\mathcal U=\prod_{i=1}^k\mathcal U^i$ where each $\mathcal  U^i$ is a set on the domain of $\xi^i$. The set of constraints
$$g(x;\xi^i)\in\mathcal A^i,i=1,\ldots,k\ \ \forall\ \xi\in\mathcal U$$
is equivalent to
$$g(x;\xi^i)\in\mathcal A^i\ \ \forall\ \xi^i\in\mathcal U^i,\ \ i=1,\ldots,k$$\label{lemma:intersection}
\end{lemma}

 Note that in approximating a joint CCP, all the $\mathcal U^i$ in Lemma \ref{lemma:intersection} need to be jointly calibrated statistically to account for the simultaneous estimation error (which can be conducted by introducing a max operation for the intersection of sets). Intuitively, with weakly correlated data across the constraints, it fares better to use a separate $\mathcal U^i$ to represent the uncertainty of each constraint rather than using a single $\mathcal U$ and projecting it. Appendix \ref{sec:comparison} provides a formal statement to support this intuition, by arguing a lower level of conservativeness in using individual ellipsoids rather than a single aggregated block-diagonal ellipsoid.

In addition, we can borrow the following statistical tools to more tightly trace an HPR, i.e., a smaller-volume prediction set:
\begin{enumerate}
\item  When data appears in multi-modal form, we can use clustering. Label the data into different clusters (using $k$-means, Gaussian mixture models, or any other techniques), form a simple set $\mathcal U_i$ like a ball or an ellipsoid for each cluster, and use the union $\bigcup_i\mathcal U_i$ as the final shape. 
\item If the high-dimensional data set has an intrinsic low-dimensional representation, we can use dimension reduction tools like principal component analysis. Suppose $\tilde\xi=M\xi+N$, where $M\in\mathbb R^{r\times m}$ and $N\in\mathbb R^r$, is a low-dimensional representation of a raw random vector $\xi\in\mathbb R^m$. Then we can use uncertainty set in the form\begin{equation}
\mathcal U=\{(M\xi-\mu)'\Sigma^{-1}(M\xi-\mu)\leq s\},\label{uncertainty PCA}
\end{equation}
where $\mu$ is the sample mean of $\tilde{\xi}$ and $\Sigma$ is a covariance estimate of $\tilde{\xi}$. Tractability is preserved by a straightforward use of existing RO results (see Theorem \ref{thm:PCA} in Appendix \ref{sec:RO results}).
\item In situations of unstructured data where clustering or dimension reduction techniques do not apply, one approach is to view each data point as a ``cluster"  by taking the union of balls each surrounding one data point. Intriguingly, this scheme coincides with the one studied in \cite{erdougan2006ambiguous} to approximate ambiguous CCP where the underlying distribution is within a neighborhood of some baseline measure.
\end{enumerate}

We provide further illustrations of these tools in Appendix \ref{sec:single_linear_ml}.

\section{Numerical Examples}\label{sec:numerics}

We present numerical examples to illustrate the performances of our RO approach. 
In all our examples,
\begin{enumerate}
\item We set $\epsilon=0.05$ and $\delta=0.05$.
\item For each setting, we repeat the experimental run $1,000$ times, each time generating a new independent data set.
\item We define $\hat{\epsilon}$ to be the estimated expected violation probability of the obtained solution. In other words, $\hat\epsilon=\hat E_D \left[ P_{violation} \right] $, where $\hat E_D[\cdot]$ refers to the empirical expectation taken among the $1,000$ data sets, and $ P_{violation}$ denotes the probability $P(g(\hat{x}(D);\xi)\notin \mathcal{A})$. For single linear CCPs with Gaussian distributed $\xi$, $ P_{violation}$ can be computed analytically. In other cases, $ P_{violation}$ is estimated using $10,000$ new independent realizations of $\xi$. For approaches that do not depend on data, e.g., SCA, we set $\hat{\epsilon}=P_{violation}$ directly.
\item We define $\hat{\delta}=\hat P_D (P_{violation} > \epsilon )$, where $\hat P_D(\cdot)$ refers to the empirical probability with respect to the $1,000$ data sets and $P_{violation}$ is similarly defined as for $\hat{\epsilon}$. For approaches that do not depend on data, the chance constraint is always satisfied and therefore we have $\hat{\delta}=0$.
	\item We denote ``Obj. Val." as the average optimal objective value of the 1,000 solutions generated from the independent data sets.
    \item When the reconstruction technique described in Section \ref{sec:enhancement} is applied, the initial guessed solution is obtained from an uncertainty set with size calibrated to be just enough to cover $(1-\epsilon)$ of the Phase 1 data.
\end{enumerate}

Recall that $d$ is the decision space dimension, $n$ is the total sample size, and $n_1$ and $n_2$ are the sample sizes for Phases 1 and 2. These numbers differ across the examples for illustration purpose. 

Moreover, we compare our RO approaches with several methods: 
\begin{enumerate}
\item Scenario approaches, including the classical SG (\cite{campi2008exact}) described in the introduction and its variant FAST (\cite{care2014fast}). FAST was introduced to reduce the sample size requirement of the classical SG. It consists of two steps, each step using $n_1$ and $n_2$ samples respectively (the notations are unified with our method for easy comparisons). The first step of FAST is similar to SG, which solves a sampled program with $n_1$ constraints and obtains a tentative solution. The second step is a detuning step to adjust the tentative solution with the help of a ``robust feasible solution", i.e., a solution feasible for any possible $\xi$. The adjusted solution is a convex combination of the tentative solution and the robust feasible solution so that the final solution satisfies the other $n_2$ sampled constraints. In our comparison, we use the minimum required sample sizes in the detuning step suggested in \cite{care2014fast} so that the total required sample size is precisely the given overall size. We compare with FAST here since the latter elicits a small sample size requirement with the help of a validation-type scheme that is similar to our approaches applied to the RO setting. 
\item DRO with first and second moment information, where the moments lie in an ellipsoidal joint confidence region. First, supposing we are given exact first and second moments, we can reformulate a distributionally robust linear chance constraint into a quadratic constraint suggested in~\cite{ghaoui2003worst}. On the other hand, using the delta method suggested in  \cite{marandi2017extending}, we can construct ellipsoidal confidence regions for the vectorized mean and covariance matrix. Combining the quadratic constraint in~\cite{ghaoui2003worst} and the ellipsoidal set in \cite{marandi2017extending}, we can use Theorem 1 (II) and Example 4 in \cite{marandi2017extending} to reformulate the DRO with  ellipsoidal moment set into a semidefinite program. We provide further details of this reformulation in Appendix \ref{sec:quadratic DRO}.

\item DRO with uncertainty set defined by a neighborhood surrounding a reference distribution measured by a $\phi$-divergence. We use the reformulation in \cite{jiang2016data} that transforms such a distributionally robust chance constraint into an ordinary chance constraint, under the reference distribution, with an adjusted tolerance level $\epsilon^*$, which then allows us to resort to SG or SAA using Monte Carlo samples (as we will see momentarily, whichever method to resort to does not quite matter in our experiments). We use the Kullback-Leibler (KL) divergence, and construct the reference distribution using kernel density estimation (with Gaussian kernel). We set the size of the KL-divergence ball by estimating the divergence using the $k$-NN estimator, a provably consistent estimator proposed in \cite{wang2009divergence,poczos2012nonparametric} (other related estimators and theoretical results are in \cite{moon2014multivariate,liu2012exponential,pal2010estimation,poczos2012nonparametric1}). We use $k=1$ in our experiments, as the experimental results indicate that the bias increases significantly as $k$ increases. Moreover, to estimate the divergence properly, we split the data into two portions $n_1$ and $n_2$, first portion used to construct the reference kernel density, second portion used for the $k$-NN divergence estimation. The reason of this split is that, otherwise, the estimation of the reference distribution and the divergence would depend on and interfere with each others, leading to estimation accuracy so poor that the divergence estimate becomes negative all the time. We provide further implementation details in Appendix \ref{sec:phi_investigation}.

\item SCA. We will state the underlying a priori distributional assumptions in using the considered SCA, which differ case-by-case.
\end{enumerate}

When applying moment-based DRO and SCA to joint CCPs, we use the Bonferroni correction (more details in the relevant examples). We also make two additional remarks. First, when comparing the objective values from different methods, since one can always translate or scale the problem by adding/multiplying constants to distort the apparent magnitudes, we mostly focus our comparisons on the direction (bigger or smaller), which is invariant under the above distortions. Second, even though we only report the point estimates of the mean objective values and $\epsilon$, $\delta$, our conclusions in comparing the objective values and constraint violation probabilities remain unchanged even if we consider the $95\%$ confidence intervals of these estimates (from the $1,000$ experimental repetitions), and we do not report the confidence intervals for the sake of succinctness. Finally, our codes are available at https://github.com/zhyhuang/Learningbased-RO.

\subsection{Test Case 1: Multivariate Gaussian on a Single Chance Constraint} \label{sec:test_case_1}
We consider a single linear CCP
\begin{equation}
	\text{minimize\ }c'x\text{\ \ subject to\ \ }P(\xi'x\leq b)\geq1-\epsilon
	\label{eq:single_CCP_ex}
\end{equation}
where $x\in \R^d$ is the decision vector, and $c \in \R^d$, $b \in \R$ are arbitrarily chosen constants. The random vector $\xi \in \R^d$ is drawn from a multivariate Gaussian distribution with an arbitrary mean (here we set it to $-c$) and an arbitrarily chosen positive definite covariance matrix. Since \eqref{eq:single_CCP_ex} is exactly solvable when the Gaussian distribution is known, we can verify that it has a bounded optimal solution.




\begin{table}[t]
	\centering
	\caption{Optimality and feasibility performances on a single $d=11$ dimensional linear CCP with Gaussian distribution for several methods, using sample size $n=120$. The true optimal value is -1196.7.} 
			
		
		
        
        \begin{tabular}{|l||l|l|l|l|l|l||l|l|l|l|l||l|l|}
\hline
                  & RO      & Recon     & SG        & FAST                 & DRO Mo & DRO KL     & SCA      \\ \hline
$n$               & 120      & 120            & 120                 & 120    & 120 & 120          & -        \\ \hline
$n_1$             & 60       & 60       & -         & 61                    & -   & 60          & -        \\ \hline
$n_2$             & 60       & 60       & -         & 59                     & -  & 60         & -        \\ \hline
Obj. Val.         & -1189.31 & -1194.87 & -1196.60  & -1193.53            &-1187.35& 0 & -1195.07 \\ \hline
$\hat{\epsilon}$  & $1.34\times 10^{-5}$ & 0.0164 & 0.090   & 0.0164    &$2.55\times 10^{-8}$ & 0  & 0.0072   \\ \hline
$\hat{\delta}$    & 0        & 0.048    & 0.957     & 0.043           & 0        & 0        & 0        \\ \hline
\end{tabular}
	
	\label{table.n=11.less}
\end{table}

			
		
		
	

\begin{table}[t]
	\centering
	\caption{Optimality and feasibility performances on a single $d=100$ dimensional linear CCP with Gaussian distribution for several methods, using sample size $n=120$. The true optimal value is -1195.3. Results on moment-based DRO are based on 30 replications due to high computational demand.} 
			
		
		
        
        \begin{tabular}{|l||l|l|l|l|l|l||l|l|l|l|l||l|l|}
\hline
                         & RO      & Recon      & SG        & FAST     & DRO Mo & DRO KL   & SCA      \\ \hline
$n$                    & 120      & 120           & 120 	     & 120    & 120     & 120             & -      \\ \hline
$n_1$                    & 60       & 60        & -        & 61        & -      & 60        & -      \\ \hline
$n_2$   			    & 60       & 60          & -   	     & 59     	& -     & 60          & -      \\ \hline
Obj. Val.			      & -832.12 & -1112.11  & unbounded& unbounded  &-1193.21& 0 & -1193.0  \\ \hline
$\hat{\epsilon}$	      & 0        & 0.0158   & -         & -     	   &0.195 & 0   & 0.0072 \\ \hline
$\hat{\delta}$  	      & 0        & 0.041    & -          & -     	       &1 & 0        & 0      \\ \hline
\end{tabular}	
	
	\label{table.n=100.less}
\end{table}

\begin{table}[t]
	\centering
	\caption{Optimality and feasibility performances on a single $d=11$ dimensional linear CCP with Gaussian distribution for several methods, using sample size $n=336$. The true optimal value is -1196.7.}
        
        \begin{tabular}{|l||l|l|l|l|l|l||l|l|l|l|l||l|l|}
\hline
                    & RO                   & Recon    & SG      & FAST           & DRO Mo & DRO KL               & SCA      \\ \hline
$n$                 & 336                  & 336      & 336     & 336         & 336      & 336          & -        \\ \hline
$n_1$               & 212                  & 212      & -       & 318         & -          & 168              & -        \\ \hline
$n_2$               & 124                  & 124      & -        & 18         & -           & 168             & -        \\ \hline
Obj. Val.           & -1190.33             & -1195.82 & -1195.67& -1195.14              &-1188.48& 0& -1195.07 \\ \hline
$\hat{\epsilon}$    & $3.47\times 10^{-6}$ & 0.0247     & 0.0331  & 0.0259   &$2.19\times 10^{-8}$ &0  & 0.0072   \\ \hline
$\hat{\delta}$      & 0                    & 0.04      & 0.056   & 0.043      & 0            &0             & 0        \\ \hline
\end{tabular}
	
	\label{table.n=11.more}
\end{table}

\begin{table}[t]
	\centering
	\caption{Optimality and feasibility performances on a single $d=100$ dimensional linear CCP with Gaussian distribution for several methods, using sample size $n=2331$. The true optimal value is -1195.3. Results on moment-based DRO are based on 30 replications due to high computational demand.}
			
		
		
        
        \begin{tabular}{|l||l|l|l|l|l|l||l|l|l|l|l||l|l|}
\hline
                   & RO      & Recon         & SG        & FAST      & DRO Mo  & DRO KL     & SCA      \\ \hline
$n$                 & 2331      & 2331      & 2331      & 2331      & 2331            & 2331            & -      \\ \hline
$n_1$              & 1318      & 1318       & -         & 2326      & -             & 1166          & -      \\ \hline
$n_2$              & 1013      & 1013       & -         & 5         & -              & 1165           & -      \\ \hline
Obj. Val.          & -1168.35 & -1194.76    & -1194.13  & -1193.85      & -1175.48  & 0        & -1193.0  \\ \hline
$\hat{\epsilon}$   & 0        & 0.0395        & 0.0428    & 0.0386 & $8.76 \times 10^{-14}$ & 0    & 0.0072 \\ \hline
$\hat{\delta}$     & 0        & 0.051         & 0.039    & 0.052             & 0       & 0           & 0      \\ \hline
\end{tabular}	
	
	\label{table.n=100.more}
\end{table}


We consider $d=11$ and $100$ as the dimension of the decision vector. Tables \ref{table.n=11.less} and \ref{table.n=100.less} show these two cases with a small sample size $n=120$, whereas Tables   \ref{table.n=11.more} and \ref{table.n=100.more} show these cases with a bigger sample size ($336$ and $2331$ respectively) so that the classical SG provides provable feasibility guarantees. In each table, we show the results for our RO using ellipsoidal uncertainty set (``RO"), our reconstructed RO (``Recon"), SG (``SG"), FAST (``FAST"), DRO with ellipsoidal moment set (``DRO Mo"), DRO with KL-divergence set (``DRO KL") and SCA (``SCA"). The last approach does not need the data and instead assumes partial a priori distributional information. 


For our RO approaches, we use ellipsoidal uncertainty sets with estimated covariance matrix for the case $d=11$ (Tables \ref{table.n=11.less} and \ref{table.n=11.more}), and diagonalized ellipsoidal sets (i.e., only using variance estimates) for $d=100$ (Tables \ref{table.n=100.less} and \ref{table.n=100.more}) to stabilize our estimates because $n_1$ is smaller than $d$ in the latter case. The tables show that the solutions from our plain RO tend to be conservative, as $\hat{\delta}=0$. Nonetheless, the reconstructed RO is less conservative across all settings, reflected by the better average optimal values and $\hat\delta$ close to the target confidence level $0.05$. In all cases, both the plain RO and the reconstructed RO give valid (i.e., confidently feasible) solutions. 

We compare our ROs with scenario approaches. When the sample size is small (Tables \ref{table.n=11.less} and \ref{table.n=100.less}), SG cannot obtain a valid solution. In the case $d=11$, it gives $\hat\delta$ much greater than $0.05$. Furthermore, in the case $d=100$, SG gives unbounded solutions in all $1,000$ replications, as the number of sampled constraints is very close to the decision dimension. For FAST, since $b$ is chosen to be positive, we can use the origin to be the robust feasible solution. Table \ref{table.n=11.less} shows that, when $d=11$, FAST gives confidently feasible solutions. The average optimal value from reconstructed RO (-1194.87) is (slightly) better than the value from FAST (-1193.53), while RO using ellipsoidal sets is more conservative (-1189.31). However, when $d=100$ (Table \ref{table.n=100.less}), the first-step problem of FAST is unbounded in all $1,000$ replications. 

 

When the sample size is adequate (Tables \ref{table.n=11.more} and \ref{table.n=100.more}), the values of $\hat\delta$ from SG being  less than or close to $0.05$ confirms the validity of the solutions. Note that in these cases FAST gives more conservative solutions than SG (This is a general consequence from the construction of FAST that is designed to have a smaller feasible region than SG under the same dataset). RO with ellipsoidal sets obtains more conservative solutions than SG, as shown by the zero $\hat\delta$'s and worse average objective values. By using reconstruction, however, the $\hat{\delta}$'s become very close to the desired confidence level $\delta=0.05$, and the average objective values are almost identical to (and slightly better than) those obtained from SG.

The above reveal that, when the sample size is large enough, SG can perform better than our RO using basic uncertainty sets. On the other hand, our RO can provide feasibility guarantees in small-sample situations where SG may fail. FAST is valid in small-sample situations, but is more likely to have unbounded solutions in high-dimensional problems than our RO. Thus, generally, our RO appears most useful for small sample sizes when compared with scenario approaches, a benefit postulated in the previous sections. It also appears that using reconstruction can boost our performance to a comparable level as SG (and hence also FAST) in situations where the latter is applicable in the shown examples. Note that our reconstruction by design can improve the objective performance compared to plain RO, whereas FAST is primarily used to reduce the sample size requirement and is necessarily more conservative than SG in terms of achieved objective value. Finally, we note that unbounded solutions in SG can potentially be avoided by adding artificial constraints. In this regard, we show in Appendix \ref{sec:append_single_example} the same example but with additional non-negativity constraints to illustrate the comparisons further.

Next, we compare with moment-based DRO. In low-dimensional cases with $d=11$, moment-based DRO gives solutions more conservative than RO using ellipsoidal sets, as shown by the larger objective values, i.e. -1187.35 (DRO) versus -1189.31 (RO) in the small-sample case (Table \ref{table.n=11.less}) and -1184.48 (DRO) versus -1190.33 (RO) in the large-sample case (Table \ref{table.n=11.more}). The conservativeness of moment-based DRO is also revealed in the small $\hat{\epsilon}$ and $\hat{\delta}=0$ in both cases. For high-dimensional problems with $d=100$, we present the performance of moment-based DRO with only 30 replications (instead of 1000) due to the large program size and consequently the demanding computational effort when solving the reformulated semidefinite programs (although the replication size is smaller, conclusions can still be drawn rigorously, i.e., the confidence intervals of the estimated $\hat\epsilon$ and $\hat\delta$ turn out to either lie completely under or above 0.05). In the small-sample size case (Table \ref{table.n=100.less}), moment-based DRO fails to provide feasible solutions ($\hat \delta=1$, i.e., obtained solutions violate the chance constraint in all 30 replications). This can be attributed to a poor estimation of the moment confidence region with small data and high dimension (Note that forming an ellipsoidal first-and-second-moment set for moment-based DRO requires estimating a covariance matrix of size $(3d+d^2)/2\times(3d+d^2)/2$, as it uses the estimation variances of the first and second moments that involve even higher-order moments, in contrast to a size of $d\times d$ in our ellipsoidal RO). When the sample size is larger (Table \ref{table.n=11.more}), moment-based DRO provides valid feasible solutions ($\hat\delta=0$). The average objective (-1175.48) is less conservative than our plain RO (-1168.35), but is more conservative than our reconstructed RO (-1194.76). 

The above observations show that, when the moment information is well estimated (i.e., the sample size is sufficient relative to the dimension), moment-based DRO provides solutions with similar conservative level as our RO using ellipsoidal sets. However, when the sample size is too small to get reasonable estimates for the moments, moment-based DRO can fail to obtain feasible solutions. Reconstructed RO appears to outperform moment-based DRO generally. The benefits of our RO approaches in small sample and the boosted performance of reconstructed RO compared to moment-based DRO are in line with our comparisons with scenario approaches.


DRO with estimated KL-divergence set suffers from general setbacks in the experiments. In all cases we considered, the kernel density estimator cannot provide a good enough reference distribution $f_0$, so that the size of the divergence ball is too big and subsequently results in conservative solutions. The construction of $f_0$ is poor due to the curse of dimensionality in kernel density estimation whose accuracy deteriorates exponentially with the dimension, as we have a relatively high dimension compared with the data size. 
On the other hand, the performance of DRO, which relies on using the adjusted tolerance level $\epsilon^*$, appears sensitive to the divergence ball size and demands a high accuracy in estimating $f_0$. 
Subsequently, the big divergence ball size leads to a zero $\epsilon^*$ in all replications, which in turn forces us to choose a solution $x$ that satisfies the safety condition $\xi' x \leq b$ for all $\xi \in \R^d$. The origin is then output as the only such feasible solution, and the objective is 0, which are shown in Tables \ref{table.n=11.less}, \ref{table.n=100.less}, \ref{table.n=11.more}, and \ref{table.n=100.more}. This indicates that DRO with KL divergence, calibrated using density estimator and the divergence estimation technique suggested in the literature, gives overly conservative solutions for our considered problems.

 
Lastly, we compare with SCA. Consider a perturbation model for $\xi$ given by $\xi=a_0+\sum_{i=1}^L\zeta_ia_i$ where $a_i\in\mathbb R^d$ for all $i=0,1,\ldots,L$ and $\zeta_i\in\mathbb R$ are independent Gaussian variable with mean $\mu_i$ and variance $s_i^2$, such that $\mu_i \in [\mu_i^-,\mu_i^+]$ and $s^2_i \leq \sigma^2_i$. A safe approximation of (\ref{eq:single_CCP_ex}) is in \cite{ben2009robust}:
\[
\text{minimize\ }c'x\text{\ \ subject to\ \ }(a_0'x-b)+\sum_{i=1}^{L} \max [a_i'x \mu_i^-,a_i'x \mu_i^+ ] + \sqrt{2 \log (1/\epsilon)} \sqrt{\sum_{i=1}^{L}\sigma_i^2(a_i'x)^2} \leq 0.
\]

To apply this SCA to (\ref{eq:single_CCP_ex}), we set $\zeta_i$ to be independent $N(0,1)$ variables, $a_0=\mu$ and $a_i$ to be the $i$-th column of $\Sigma^{1/2}$, and $\mu_i^-=\mu_i^+=0$ and $\sigma_i^2=1$ for $i=1,...,d$.  This in fact assumes knowledge on the mean and covariance of the Gaussian vector $\xi$, thus giving an upper hand to SCA.


Tables \ref{table.n=11.less}, \ref{table.n=11.more}, \ref{table.n=100.less} and \ref{table.n=100.more} all show that the optimal objective values obtained from SCA (-1195.07 and -1193.0 respectively for $d=11,100$) are close to the true optimal values (-1196.7 and -1195.3) compared to other methods. Our ROs using ellipsoidal sets obtain more conservative solutions  generally. The relative conservativeness also shows up in reconstructed RO with small sample sizes (Tables \ref{table.n=11.less} and \ref{table.n=100.less}), but with more samples (Tables \ref{table.n=11.more} and \ref{table.n=100.more}) our reconstructed RO outperforms the considered SCA.

Note that in this example the normality, and the mean and covariance information used in the SCA, makes the latter perform very well. Our RO using estimated ellipsoidal sets does not achieve this level of preciseness. However, the reconstructed RO can still outperform this SCA when the sample size is large enough. Note that the performance of SCA depends on the true distribution (as it is related to the tightness of the SCA constraint in approximating the chance constraint). In the next example, we consider an alternate underlying distribution where SCA does not perform as well.

\subsection{Test Case 2: Beta Models on a Single Chance Constraint}\label{sec:test_case_2} 
We consider the single linear CCP in \eqref{eq:single_CCP_ex}, where each component of $\xi$ is now bounded. We use a perturbation model for $\xi$ given by $\xi=a_0+\sum_{i=1}^L\zeta_ia_i$ where $a_i\in\mathbb R^d$ for all $i=0,1,\ldots,L$ and $\zeta_i\in\mathbb R$ are independent random variables each with mean zero and bounded in $[-1,1]$, where $d=10$, $L=10$ and $a_i \in \R^{10}$ being known arbitrarily chosen vectors. This allows the use of an SCA stated below. In particular, we set each $\zeta_i$ to be a Beta distribution with parameters $\alpha=10$ and $\beta=10$ that is multiplied by 2 and shifted by 1. Similar to Section \ref{sec:test_case_1}, we set $c$ to be the negative of the mean of $\xi$ and $b\in \R$ is an arbitrarily chosen positive constant. 

Regarding the comparison with SCA, this problem is supplementary to the Gaussian cases in Section \ref{sec:test_case_1} in that it presents performances of SCA when we use less information about $\xi$. Suppose that we have chosen a correct perturbation model in the SCA (i.e., knowledge of $d,L,a_i$ and the boundedness on $[-1,1]$). We use the Hoeffding inequality to replace the chance constraint with $\eta\sqrt{\sum_{i=1}^L(a_i'x)^2}\leq b-a_0'x$, where $\eta\geq\sqrt{2\log(1/\epsilon)}$. This SCA is equivalent to an RO imposing an uncertainty set $\mathcal U=\{\zeta:\|\zeta\|_2\leq\eta\}$ where $\zeta=(\zeta_i)_{i=1,\ldots,L}'$ is the vector of perturbation random variables (\cite{ben2009robust} Section 2.3).



\begin{table}[ht]
	\centering
    \captionsetup{width=.9\linewidth}
    
\caption{Optimality and feasibility performances on a single  $d=10$ dimensional linear CCP with the Beta-perturbation model for several methods, using sample size $n=120$.}
    \begin{tabular}{|l|l|l|l|l|l|l||l|l|}
\hline
                & RO    & Recon         & SG       & FAST               & DRO Mo & DRO KL   & SCA   \\ \hline
$n$             & 120      & 120       & 120      & 120                & 120      & 120           & -       \\ \hline
$n_1$           & 60       & 60        & -        & 61                 & -        & 60         & -       \\ \hline
$n_2$           & 60       & 60        & -        & 59                  & -     & 60        & -       \\ \hline
Obj. Val.       & -988.78 & -1087.85    & -1114.57 & -1071.77          &-968.30 &0   & -815.06 \\ \hline
$\hat\epsilon$  & $1.02\times10^{-5}$  & 0.0161 & 0.0643  & 0.0171   &0 &0         & 0       \\ \hline
$\hat{\delta}$  & 0        & 0.037     & 0.723    & 0.063                & 0    &0          & 0       \\ \hline
\end{tabular}
    
		\vspace{3mm}
			
		\label{table:bounded_pert}
        \vspace{-5mm}
	\end{table}

Table \ref{table:bounded_pert} shows the results from different approaches with sample size $n=120$. Our RO performs better than SCA in terms of achieved objective values  ($-988.78$ against $-815.06$), the latter appearing more conservative than the example in Section \ref{sec:test_case_1} as shown by $\hat\epsilon=0$. Also, as in the previous example, reconstruction boosts further our RO performance (from $-988.78$ to $-1087.85$). Our RO here performs better than SCA because the latter, derived on a worst-case basis, does not tightly apply to the ``truth" in this example, i.e., the Hoeffding bound does not lead to tight performance guarantees on the scaled Beta distribution (putting aside the assumed knowledge of $d,L,a_i$ and the boundedness on $[-1,1]$ when applying the SCA). Note that, since SCA also has an RO interpretation, the above observations show the superiority of our geometry or size selection of the uncertainty set. Our fully nonparametric approach shows full-fledged advantage than SCA in this example.




We also report the outcomes of SG, which breaks down as shown by $\hat\delta$ being much bigger than $0.05$, as 120 observations is not enough to achieve the needed feasibility confidence. FAST obtains valid solutions, and outperforms our RO with ellipsoidal sets but underperforms our reconstructed RO in terms of achieved objective value. Moment-based DRO also obtains valid solutions, but is conservative as shown by $\hat\delta=0$ and $\hat\epsilon=0$. Its objective value underperforms our RO approaches. For divergence-based DRO, the poor construction of a reference distribution again leads to a large divergence ball size, which renders the adjusted tolerance level $\epsilon^{*}$ to be 0 in all but one out of 1000 replications (for the one replication where $\epsilon^*$ is non-zero, it is $\epsilon^{*}=1.10\times 10^{-11}$) and essentially outputs the origin as the solution all the time. In this example, our reconstructed RO performs the best among all considered approaches. 


\subsection{Test Case 3: Multivariate Gaussian on Joint Chance Constraints} 
\label{sec:test_case_3}


We consider a joint CCP with $d=11$ variables and $l=15$ constraints in the form 
\begin{equation}
	\text{minimize\ }c'x\text{\ \ subject to\ \ }P(A x \leq b) \geq 1-\epsilon,\ x\geq 0\label{joint CCP opt}
\end{equation}
where $c \in \R^{11}$ and $b\in \R^{15}$ are arbitrary constants, and $b$ is positive in each element. The random vector $\xi= vec(A)$ is generated from a multivariate Gaussian distribution with mean $vec(\bar{A})$ and covariance matrix $\Sigma$, where $\bar{A} \in \R^{15\times 11}$ is arbitrary and $\Sigma \in \R^{165 \times 165}$ is also an arbitrary positive definite matrix. 

Tables \ref{table:joint.less} and \ref{table:joint.more} present the experimental results using two different sample sizes on the same problem. We use diagonalized ellipsoids in our RO, and conduct reconstruction with scaling parameters $k_i$ described in Appendix \ref{sec:missing test case}. To use DRO and SCA, we apply the Bonferroni correction to decompose the joint CCP, by evenly dividing the tolerance level into $\epsilon/m$ to create individual chance constraints. For each individual chance constraint, we construct DRO and SCA constraint following the scheme in Section \ref{sec:test_case_1}.

\begin{table}
	\centering
	\caption{Optimality and feasibility performances on a joint linear CCP with Gaussian distribution for several methods, using sample size $n=120$.}
			
			
			
\begin{tabular}{|l||l|l|l|l|l|l||l|l|l|l|l||l|l|}
\hline
                  & RO                    & Recon    & SG       & FAST              & DRO Mo   & DRO KL      & SCA      \\ \hline
$n$               & 120                   & 120      & 120      & 120         & 120       & 120           & -        \\ \hline
$n_1$             & 60                    & 60       & -        & 61          & -          & 60        & -        \\ \hline
$n_2$             & 60                    & 60       & -        & 59          & -         & 60             & -        \\ \hline
Obj. Val.         & -6956.49              & -7920.12 & -9283.35 & -8925.74       &-3996.87   & 0  & -8927.71 \\ \hline
$\hat{\epsilon}$  & $3.46 \times 10^{-5}$ & 0.0161   & 0.0581   & 0.0169            & 0  & 0     & 0.026    \\ \hline
$\hat{\delta}$    & 0                     & 0.044    & 0.607    & 0.045             & 0     & 0       & 0        \\ \hline
\end{tabular}

	
	\label{table:joint.less}
\end{table}

\begin{table}
	\centering
	\caption{Optimality and feasibility performances on a joint linear CCP with Gaussian distribution for several methods, using sample size $n=336$.}
			
			
			
\begin{tabular}{|l||l|l|l|l|l|l||l|l|l|l|l||l|l|}
\hline
                & RO                   & Recon         & SG        & FAST           & DRO Mo & DRO KL         & SCA      \\ \hline
$n$             & 336                   & 336          & 336       & 336         & 336      & 336      & -        \\ \hline
$n_1$           & 212                   & 212          & -         & 318         & -          &168        & -        \\ \hline
$n_2$           & 124                   & 124          & -         & 18          & -         &168         & -        \\ \hline
Obj. Val.       & -7146.54              & -8029.83     & -9130.95  & -9081.81     &-4209.86  & 0   & -8927.71 \\ \hline
$\hat{\epsilon}$& $7.32 \times 10^{-5}$ & 0.0235       & 0.0223    & 0.0185    & 0           & 0     & 0.026    \\ \hline
$\hat{\delta}$  & 0                     & 0.038        & 0.005     & 0.002       & 0          & 0        & 0        \\ \hline
\end{tabular}

	
	\label{table:joint.more}
\end{table}

Comparing with scenario approaches, we see that, much like the examples in Sections \ref{sec:test_case_1} and \ref{sec:test_case_2}, SG fails with small sample size (confirmed by $\hat{\delta}$ much larger than 0.05 in Table \ref{table:joint.less}), but obtains valid solutions as sample size grows (confirmed by $\hat\delta<0.05$ in Table \ref{table:joint.more}).  While reconstruction improves the optimal values for RO in both cases, SG (and so is FAST) gives better optimal value ($-9130.95$) than reconstructed RO ($-8029.83$) under a big sample size. Moment-based DRO appears very conservative for both small and large sample cases, as the obtained average objective values (-3996.87 and -4209.86) are much greater than other approaches, including our ROs, and the associated $\hat\epsilon$ and $\hat\delta$ are 0. Like the previous experiments, divergence-based DRO outputs the origin as the solution and gives objective value 0 due to over-sized uncertainty sets. On the other hand, SCA obtains a better solution than our ROs, thanks to the tightness of the approximation for Gaussian distributions.





\subsection{Test Case 4: Beta Models on Joint Chance Constraints} 

We consider the joint CCP in \eqref{joint CCP opt} with a bounded random vector $\xi$. We use the perturbation model described in Section \ref{sec:test_case_2}, where $d=165$, $L=165$ and $a_i \in \R^{165}, i=1,...,L$ are arbitrarily chosen vectors, and  the same random variables for $\zeta_i$'s as in Section \ref{sec:test_case_2}. Again, we apply the Bonferroni correction to invoke DRO and SCA as in Section \ref{sec:test_case_3}, and the corresponding schemes for each individualized chance constraint as in Section \ref{sec:test_case_2}.

\begin{table}
	\centering
	\caption{Optimality and feasibility performances on a joint linear CCP with beta distribution for several methods, using sample size $n=120$.}
			
			
			
\begin{tabular}{|l||l|l|l|l|l|l||l|l|l|l|l||l|l|}
\hline
                  & RO                    & Recon     & SG         & FAST        & DRO Mo  & DRO KL         & SCA      \\ \hline
$n$              & 120                    & 120      & 120        & 120         & 120   & 120              & -        \\ \hline
$n_1$            & 60                     & 60       & -          & 61          & -        & 60           & -        \\ \hline
$n_2$             & 60                     & 60       & -          & 59          & -       & 60           & -        \\ \hline
Obj. Val.         & -1241.05               & -1796.74 & -2105.77   & -1732.73      &-230.74 &0  & -361.079 \\ \hline
$\hat{\epsilon}$  & $6.96 \times 10^{-5}$ & 0.0138    & 0.0577     & 0.0170      & 0       & 0          & 0        \\ \hline
$\hat{\delta}$    & 0                      & 0.022    & 0.576      & 0.045       & 0      & 0          & 0       \\ \hline
\end{tabular}

	
	\label{table:joint.beta.less}
\end{table}

\begin{table}
	\centering
	\caption{Optimality and feasibility performances on a joint linear CCP with beta distribution for several methods, using sample size $n=336$.}

\begin{tabular}{|l||l|l|l|l|l|l||l|l|l|l|l||l|l|}
\hline
                     & RO                    & Recon           & SG          & FAST       & DRO Mo   & DRO KL     & SCA      \\ \hline
$n$                 & 336                    & 336            & 336          & 336       & 336   & 336      & -        \\ \hline
$n_1$                & 212                    & 212            & -            & 318        & - &168             & -        \\ \hline
$n_2$               & 124                    & 124            & -            & 18           & -    &168        & -        \\ \hline
Obj. Val.           & -1304.89               & -1911.36       & -1881.69     & -1828.98   &-251.69& 0 & -361.079 \\ \hline
$\hat{\epsilon}$    & $1.20 \times 10^{-4}$ & 0.0199        & 0.0229      & 0.0192    & 0      & 0        & 0        \\ \hline
$\hat{\delta}$      & 0                      & 0.023        & 0.004          & 0.003     & 0      & 0        & 0      \\ \hline
\end{tabular}

	
	\label{table:joint.beta.more}
\end{table}
Tables \ref{table:joint.beta.less} and \ref{table:joint.beta.more} show our experimental results. The major difference with Section \ref{sec:test_case_3} is that now our reconstructed RO outperforms all other methods including SG and SCA: It gives smaller objective values than FAST under both small and big sample sizes. It also gives smaller objective values than SG under big sample size, while SG does not give valid solutions under small sample size. SCA is very conservative in this case, and DROs (both moment- and divergence-based) continue to be very conservative, all of whom our RO significantly outperforms.

\subsection{Test Case 5: $t$- and Log-Normal Distributions}

We consider problems with two heavier-tailed distributions, namely $t$- and lognormal. We test both the single CCP \eqref{eq:single_CCP_ex} and the joint CCP \eqref{joint CCP opt} with different dimensions and sample sizes. Since the considered SCA does not apply to these distributions, we do not include it in our comparisons here.

Tables \ref{table.T.n=11.less}, \ref{table.T.n=11.more} and \ref{table.single.T.n=100} show the comparisons among different approaches for the single CCP, and Tables \ref{table.joint.T.n=11.less} and \ref{table.joint.T.n=11.more} show the counterparts for joint CCP, when $\xi$ is generated from a multivariate $t$-distribution with degree of freedom 5 and an arbitrary positive definite dispersion matrix. The comparisons are largely consistent with the Gaussian and beta cases shown in the previous subsections. Compared with SG, our ROs output feasible solutions in the small-sample case ($n=120$), whereas SG struggles to obtain feasible solutions ($\hat\delta$ much greater than 0.05 in Tables \ref{table.T.n=11.less} and \ref{table.joint.T.n=11.less}). In the large-sample case ($n=336$), SG gains enough feasibility and outperforms our plain RO in average objective value (-1175.04 versus -1126.66 in the single CCP case in Table \ref{table.T.n=11.more}, and -7387.98 versus -5778.44 in the joint CCP case in Table \ref{table.joint.T.n=11.more}), but underperforms our reconstructed RO (-1175.64 and -7562.60 for single and joint CCPs respectively). FAST remedies the infeasibility issue of SG in the small-sample cases and outperforms our plain RO. On the other hand, our reconstructed RO performs competitively against FAST. Among all four cases where $d=11$, reconstructed RO outperforms FAST in three cases but underperforms in the case of small-sample joint CCP (average objective values -1166.52, -1175.64 and -7562.60 versus -1158.27, -1170.35 and -7173.97 in  Tables \ref{table.T.n=11.less}, \ref{table.T.n=11.more} and \ref{table.joint.T.n=11.more} respectively, and -6499.93 versus -7220.37 in Table \ref{table.joint.T.n=11.less}). Note that, when the dimension is large ($d=100$ in Table \ref{table.single.T.n=100}), SG and FAST output unbounded solutions in all 1000 experimental replications, whereas plain and reconstructed RO output feasible bounded solutions.

Like in the previous subsections, our reconstructed RO outperforms moment-based DRO in all cases. When the dimension is large ($d=100$ in Table \ref{table.single.T.n=100}), moment-based DRO fails to obtain feasible solutions in all 30 replications, attributed to the difficulty in estimating valid moment confidence regions. Compared to our plain RO, moment-based DRO outperforms in single CCP (-1134.38 and -1137.19 versus -1112.75 and -1126.66 in Tables \ref{table.T.n=11.less} and \ref{table.T.n=11.more} respectively), but underperforms in joint CCP (-3888.63 and -3891.83 versus -4229.6 and -5778.44 in Tables \ref{table.joint.T.n=11.less} and \ref{table.joint.T.n=11.more} respectively). Lastly, divergence-based DRO is once again very conservative, resulting in zero objective values all the time. 




\begin{table}[t]
	\centering
	\caption{Optimality and feasibility performances on a single $d=11$ dimensional linear CCP with $t$-distribution for several methods, using sample size $n=120$.}

        \begin{tabular}{|l||l|l|l|l|l|l|l||l|l|l|l|l||l|l|}
\hline
                  & RO      & Recon     & SG        & FAST                & DRO Mo & DRO KL          \\ \hline
$n$               & 120      & 120      & 120       & 120                 & 120   & 120              \\ \hline
$n_1$             & 60       & 60       & -         & 61                  & -     & 60               \\ \hline
$n_2$             & 60       & 60       & -         & 59                  & -     & 60              \\ \hline
Obj. Val.         & -1112.75 & -1166.52  & -1182.20  & -1158.27        &-1134.38 & 0  \\ \hline
$\hat{\epsilon}$  & 0.000252 & 0.0161 & 0.0910   & 0.0172    &0.000461& 0    \\ \hline
$\hat{\delta}$    & 0        & 0.046    & 0.961     & 0.064            & 0    & 0              \\ \hline
\end{tabular}
	
	\label{table.T.n=11.less}
\end{table}

\begin{table}[t]
	\centering
	\caption{Optimality and feasibility performances on a single $d=11$ dimensional linear CCP with $t$-distribution for several methods, using sample size $n=336$.}

        \begin{tabular}{|l||l|l|l|l|l|l|l||l|l|l|l|l||l|l|}
\hline
                  & RO      & Recon     & SG        & FAST                & DRO Mo & DRO KL            \\ \hline
$n$               & 336      & 336      & 336       & 336                 & 336   & 336                  \\ \hline
$n_1$             & 212       & 212       & -         & 318                  & -   & 168                 \\ \hline
$n_2$             & 124       & 124       & -         & 18                  & -     & 168                 \\ \hline
Obj. Val.         & -1126.66 & -1175.64 & -1175.04  & -1170.35          &-1137.19 & 0   \\ \hline
$\hat{\epsilon}$  & 0.00023 & 0.024 & 0.0334  & 0.0259  &0.000407 & 0   \\ \hline
$\hat{\delta}$    & 0        & 0.055    & 0.069    & 0.04               & 0    & 0              \\ \hline
\end{tabular}
	
	\label{table.T.n=11.more}
\end{table}	

\begin{table}[t]
	\centering
	\caption{Optimality and feasibility performances on a single $d=100$ dimensional linear CCP with $t$-distribution for several methods, using sample size $n=120$. Results on moment-based DRO are based on 30 replications due to high computational demand.}

        \begin{tabular}{|l||l|l|l|l|l|l|l||l|l|l|l|l||l|l|}
\hline
                  & RO      & Recon     & SG        & FAST                & DRO Mo & DRO KL           \\ \hline
$n$               & 120      & 120      & 120       & 120                 & 120   & 120                \\ \hline
$n_1$             & 60       & 60       & -         & 61                  & -      & 60             \\ \hline
$n_2$             & 60       & 60       & -         & 59                  & -    & 60           \\ \hline
Obj. Val.         & -1077.56 & -1184.45 & unbounded  & unbounded           & -1190.70 & 0\\ \hline
$\hat{\epsilon}$  & $6.00\times 10^{-7}$ & 0.0156 & -  & -  & 0.22 & 0     \\ \hline
$\hat{\delta}$    & 0        & 0.045    & -     & -               & 1   & 0           \\ \hline
\end{tabular}
	
	\label{table.single.T.n=100}
\end{table}	

\begin{table}[t]
	\centering
	\caption{Optimality and feasibility performances on a joint $d=11$ dimensional linear CCP with $t$-distribution for several methods, using sample size $n=120$.}

        \begin{tabular}{|l||l|l|l|l|l|l|l||l|l|l|l|l||l|l|}
\hline
                  & RO      & Recon     & SG        & FAST                & DRO Mo & DRO KL         \\ \hline
$n$               & 120      & 120      & 120       & 120                 & 120   & 120                \\ \hline
$n_1$             & 60       & 60       & -         & 61                  & -      & 60             \\ \hline
$n_2$             & 60       & 60       & -         & 59                  & -    & 60           \\ \hline
Obj. Val.         & -4229.6  & -6499.93  & -8313  & -7220.37           &-3888.63&0  \\ \hline
$\hat{\epsilon}$  & 0.00108 & 0.00847 & 0.0404  & 0.0152   &$4.17\times 10^{-4}$ &0  \\ \hline
$\hat{\delta}$    & 0        & 0.002    & 0.284     & 0.048               & 0   & 0           \\ \hline
\end{tabular}
	
	\label{table.joint.T.n=11.less}
\end{table}

\begin{table}[t]
	\centering
	\caption{Optimality and feasibility performances on a joint $d=11$ dimensional linear CCP with $t$-distribution for several methods, using sample size $n=336$.}

        \begin{tabular}{|l||l|l|l|l|l|l|l||l|l|l|l|l||l|l|}
\hline
                  & RO      & Recon     & SG        & FAST                & DRO Mo & DRO KL           \\ \hline
$n$               & 336      & 336      & 336       & 336                 & 336   & 336           \\ \hline
$n_1$             & 212       & 212       & -         & 318                  & -    & 168              \\ \hline
$n_2$             & 124       & 124       & -         & 18                  & -     & 168              \\ \hline
Obj. Val.         & -5778.44 & -7562.60 & -7387.98  & -7173.97           & -3891.83& 0\\ \hline
$\hat{\epsilon}$  & 0.00248 & 0.0133 & 0.0144  & 0.0126  & $3.97\times 10^{-4}$ & 0     \\ \hline
$\hat{\delta}$    & 0        & 0    & 0     & 0               & 0   & 0           \\ \hline
\end{tabular}
	
	\label{table.joint.T.n=11.more}
\end{table}


Next we consider $\xi$ generated from log-normal distributions with arbitrarily chosen means and covariance matrices. Tables \ref{table.lognormal.n=11.less}, \ref{table.lognormal.n=11.more} and \ref{table.single.lognormal.n=100} show the results for the single CCP, while Tables \ref{table.joint.lognormal.n=11.less} and \ref{table.joint.lognormal.n=11.more} show those for the joint CCP. The comparisons are quite similar to the $t$-distribution cases. SG in small sample outputs invalid solutions ($\hat \delta$ much greater than 0.05), and in large sample outputs solutions with average objective values (e.g. -683.60 in Table \ref{table.lognormal.n=11.more}) better than our plain RO (-354.10) but worse than our reconstructed RO (-685.01). FAST remedies the infeasibility issue of SG in the small-sample cases, but underperforms our reconstructed RO in all cases. Moment-based DRO outperforms our plain RO but underperforms our reconstructed RO in all cases, and it continues to struggle in obtaining feasible solutions for high-dimensional problems ($\hat{\delta}=1$ in Table \ref{table.single.lognormal.n=100}). Lastly, divergence-based DRO continues to be conservative and outputs zero objective values. In all considered settings, reconstructed RO appears the best among all compared methods in terms of feasibility and optimality.

\begin{table}[t]
	\centering
	\caption{Optimality and feasibility performances on a single $d=11$ dimensional linear CCP with log-normal distribution for several methods, using sample size $n=120$.}

        \begin{tabular}{|l||l|l|l|l|l|l|l||l|l|l|l|l||l|l|}
\hline
                  & RO      & Recon     & SG        & FAST                & DRO Mo & DRO KL      \\ \hline
$n$               & 120      & 120      & 120       & 120                 & 120  & 120             \\ \hline
$n_1$             & 60       & 60       & -         & 61                  & -     & 60              \\ \hline
$n_2$             & 60       & 60       & -         & 59                  & -      & 60           \\ \hline
Obj. Val.         & -294.00  & -588.58  & -784.27  & -510.38         &-418.30&0  \\ \hline
$\hat{\epsilon}$  & $1.45\times 10^{-4}$ & 0.0164 & 0.0902   & 0.0159  &$5.11\times 10^{-4}$&0   \\ \hline
$\hat{\delta}$    & 0        & 0.041    & 0.961     & 0.048               & 0   &  0            \\ \hline
\end{tabular}
	
	\label{table.lognormal.n=11.less}
\end{table}

\begin{table}[t]
	\centering
	\caption{Optimality and feasibility performances on a single $d=11$ dimensional linear CCP with log-normal distribution for several methods, using sample size $n=336$.}

        \begin{tabular}{|l||l|l|l|l|l|l|l||l|l|l|l|l||l|l|}
\hline
                  & RO      & Recon     & SG        & FAST                & DRO Mo & DRO KL          \\ \hline
$n$               & 336      & 336      & 336       & 336                 & 336   & 336                \\ \hline
$n_1$             & 212       & 212       & -         & 318                  & -     &168                \\ \hline
$n_2$             & 124       & 124       & -         & 18                  & -     &168                \\ \hline
Obj. Val.         & -354.10 & -685.01 & -683.60  & -646.83         &-429.75&0  \\ \hline
$\hat{\epsilon}$  & $8.07\times 10^{-5}$& 0.0243 & 0.0333 & 0.0261   & $3.33\times 10^{-4}$&0  \\ \hline
$\hat{\delta}$    & 0        & 0.057    & 0.052    & 0.033               & 0     &0           \\ \hline
\end{tabular}
	
	\label{table.lognormal.n=11.more}
\end{table}	

\begin{table}[t]
	\centering
	\caption{Optimality and feasibility performances on a single $d=100$ dimensional linear CCP with log-normal distribution for several methods, using sample size $n=120$. Results on moment-based DRO are based on 30 replications due to high computational demand.}

        \begin{tabular}{|l||l|l|l|l|l|l|l||l|l|l|l|l||l|l|}
\hline
                  & RO      & Recon     & SG        & FAST                & DRO Mo & DRO KL           \\ \hline
$n$               & 120      & 120      & 120       & 120                 & 120   & 120                \\ \hline
$n_1$             & 60       & 60       & -         & 61                  & -      & 60             \\ \hline
$n_2$             & 60       & 60       & -         & 59                  & -    & 60           \\ \hline
Obj. Val.         & -309.93 & -784.24 & unbounded  & unbounded           & -1030.52 & 0\\ \hline
$\hat{\epsilon}$  & $6.00\times 10^{-6}$ & 0.0174 & -  & -  & 0.2772 & 0     \\ \hline
$\hat{\delta}$    & 0        & 0.063    & -     & -               & 1   & 0           \\ \hline
\end{tabular}
	
	\label{table.single.lognormal.n=100}
\end{table}

\begin{table}[t]
	\centering
	\caption{Optimality and feasibility performances on a joint $d=11$ dimensional linear CCP with log-normal distribution for several methods, using sample size $n=120$.}

        \begin{tabular}{|l||l|l|l|l|l|l|l||l|l|l|l|l||l|l|}
\hline
                  & RO      & Recon     & SG        & FAST                & DRO Mo& DRO KL        \\ \hline
$n$               & 120      & 120      & 120       & 120                 & 120   & 120              \\ \hline
$n_1$             & 60       & 60       & -         & 61                  & -      & 60              \\ \hline
$n_2$             & 60       & 60       & -         & 59                  & -       & 60           \\ \hline
Obj. Val.         & -0.1284  & -1.1166  & -4.5359  & -1.0369          & -0.8360& 0\\ \hline
$\hat{\epsilon}$  & 0.00228 & 0.0157& 0.0598  & 0.0165   & 0.0131& 0  \\ \hline
$\hat{\delta}$    & 0        & 0.043    & 0.646    & 0.044           &  0.006 & 0          \\ \hline
\end{tabular}
	
	\label{table.joint.lognormal.n=11.less}
\end{table}

\begin{table}[t]
	\centering
	\caption{Optimality and feasibility performances on a joint $d=11$ dimensional linear CCP with log-normal distribution for several methods, using sample size $n=336$.}

        \begin{tabular}{|l||l|l|l|l|l|l|l||l|l|l|l|l||l|l|}
\hline
                  & RO      & Recon     & SG        & FAST                & DRO Mo& DRO KL          \\ \hline
$n$               & 336      & 336      & 336       & 336                 & 336    & 336               \\ \hline
$n_1$             & 212       & 212       & -         & 318                  & -     &168                \\ \hline
$n_2$             & 124       & 124       & -         & 18                  & -     &168                \\ \hline
Obj. Val.         & -0.0844 & -1.9373 & -1.7135  & -1.4058     &-1.2021&0  \\ \hline
$\hat{\epsilon}$  & 0.0074 & 0.0239 & 0.0238 & 0.0197 &0.0131&0    \\ \hline
$\hat{\delta}$    & 0        & 0.05    & 0.011     & 0.007        & 0.026 &0            \\ \hline
\end{tabular}
	
	\label{table.joint.lognormal.n=11.more}
\end{table}

\subsection{Summary on the Experiment Results}
From the results in this section (and additional ones in Appendix \ref{sec:numerics added}), we highlight the following situations where our method is the most recommended. 

The competitiveness of our method compared with scenario approaches is most seen in small-sample situations. Classical SG needs a much larger sample size than ours to achieve feasibility. FAST is capable of obtaining feasible solutions in small-sample cases, but appears more susceptible than RO in generating unbounded solutions. With reconstruction, our approach tends to work as well as SG and FAST for large sample (when they are all applicable). Moreover, our reconstruction has the capability to improve the optimality over plain RO, whereas FAST is by design always more conservative than SG in terms of optimality. Nonetheless, we should mention that some constraint removal approaches like sampling-and-discarding (\cite{campi2011sampling}) can improve SG performances in large-sample situations.

Compared to our ROs, moment-based DRO can generate infeasible solutions when the problem dimension is high compared to data size (e.g., $d=100$ and $n=120$), attributed to the difficulty in constructing valid moment confidence regions. In cases where moment-based DRO generates valid solutions, the solution performances seem to be sometimes better, sometimes worse than our plain RO, but in all considered instances they perform worse than our reconstructed RO. KL-divergence-based DRO appears to perform poorly in the experiments due to the challenge in obtaining a small enough divergence ball size (To get a further sense of this behavior, we investigate a very low-dimensional problem ($d=3$) with sufficient sample size in Section \ref{sec:phi_investigation}, where divergence-based DRO provides nontrivial but still conservative solutions).


Lastly, compared with SCA, our performance is best seen when the data is non-normal. In this case the approximate constraint in SCA may not tightly approximate the original chance constraint and tends to be significantly more conservative than our approach. Moreover, SCA generally requires at least some partial distributional knowledge (e.g., moments, support) in deriving the needed relaxing constraint, in contrast to our approach that is fully data-driven and nonparametric.

\ACKNOWLEDGMENT{A preliminary version of this work appeared in the Proceedings of the Winter Simulation Conference 2016. The research of L. Jeff Hong is supported in part by the Hong Kong Research Grants Council under grant GRF 11504017 and the Natural Science Foundation of China under grants No.~71991470 and No.~71991473. The research of Zhiyuan Huang and Henry Lam is supported in part by the National Science Foundation under grants CMMI-1542020, CMMI-1523453, CAREER CMMI-1834710 and IIS-1849280.}


\bibliographystyle{informs2014}
\bibliography{bibliography1}



\ECSwitch


\ECHead{Appendix}\label{sec:proof}


\section{Missing Proofs in Section \ref{sec:proposed}}\label{sec:proof}
\proof{Proof of Theorem \ref{consistency}.}
\underline{Proof of \ref{thm pt2}.}
Let $Bin(n,p)$ be a binomial variable with number of trials $n$ and success probability $p$. Then \eqref{quantile choice} can be written as
\begin{equation}
i^*=\min\left\{r:P(Bin(n_2,1-\epsilon)\leq r-1)\geq1-\delta,\ 1\leq r\leq n_2\right\}\label{quantile choice1}
\end{equation}			
Note that by the Berry-Essen Theorem,
\begin{eqnarray}
&&P(Bin(n_2,1-\epsilon)\leq r-1)-\Phi\left(\frac{r-1-n_2(1-\epsilon)}{\sqrt{n_2(1-\epsilon)\epsilon}}\right)\notag\\
&=&P\left(\frac{Bin(n_2,1-\epsilon)-n_2(1-\epsilon)}{\sqrt{n_2(1-\epsilon)\epsilon}}\leq\frac{r-1-n_2(1-\epsilon)}{\sqrt{n_2(1-\epsilon)\epsilon}}\right)-\Phi\left(\frac{r-1-n_2(1-\epsilon)}{\sqrt{n_2(1-\epsilon)\epsilon}}\right)\notag\\
&=&O\left(\frac{1}{\sqrt{n_2}}\right)\label{Berry-Essen}
\end{eqnarray}
uniformly over $r\in\mathbb N^+$, where $\Phi$ is the distribution function of standard normal. Since $i^*$ in \eqref{quantile choice1} is chosen such that $P(Bin(n_2,1-\epsilon)\leq i^*-1)\geq1-\delta$ (where we define $i^*=n_2+1$ if no choice of $r$ is valid), we have, for any $\gamma>0$, $i^*$ satisfies
$$\Phi\left(\frac{i^*-1-n_2(1-\epsilon)}{\sqrt{n_2(1-\epsilon)\epsilon}}\right)+\gamma\geq1-\delta$$
for large enough $n_2$, which gives
\begin{equation}
i^*\geq1+n_2(1-\epsilon)+\sqrt{n_2(1-\epsilon)\epsilon}\Phi^{-1}(1-\delta-\gamma)\label{interim3}
\end{equation}
for large enough $n_2$.

On the other hand, we claim that $i^*$ also satisfies, for any $\gamma>0$,
\begin{equation}
\Phi\left(\frac{i^*-1-n_2(1-\epsilon)}{\sqrt{n_2(1-\epsilon)\epsilon}}\right)\leq1-\delta+\gamma\label{interim2}
\end{equation}
for large enough $n_2$. If not, then there exists an $\gamma>0$ such that
$$\Phi\left(\frac{i^*-1-n_2(1-\epsilon)}{\sqrt{n_2(1-\epsilon)\epsilon}}\right)>1-\delta+\gamma$$
infinitely often, which implies
$$P(Bin(n_2,1-\epsilon)\leq i^*-1)+O\left(\frac{1}{\sqrt{n_2}}\right)>1-\delta+\gamma$$
or
$$P(Bin(n_2,1-\epsilon)\leq i^*-1)>1-\delta+\tilde\gamma$$
infinitely often for some $0<\tilde\gamma<\gamma$. By the choice of $i^*$, we conclude that there is no $r$ that satisfies
$$1-\delta\leq P(Bin(n_2,1-\epsilon)\leq r-1)\leq1-\delta+\tilde\gamma$$
infinitely often, which is impossible. Therefore, \eqref{interim2} holds for large enough $n_2$, and we have
\begin{equation}
i^*\leq1+n_2(1-\epsilon)+\sqrt{n_2(1-\epsilon)\epsilon}\Phi^{-1}(1-\delta+\gamma)\label{interim4}
\end{equation}

Combining \eqref{interim3} and \eqref{interim4}, and noting that $\gamma$ is arbitrary, we have
\begin{equation}
\sqrt{n_2}\left(\frac{i^*}{n_2}-(1-\epsilon)\right)\to\sqrt{(1-\epsilon)\epsilon}\Phi^{-1}(1-\delta)\label{interim6}
\end{equation}
almost surely. The same argument also shows that $i^*$ is well-defined for large enough $n_2$ almost surely.

It suffices to show that
\begin{equation}
\mathbb P_{D_2}(1-\epsilon-\gamma\leq P(\xi\in\mathcal U)\leq1-\epsilon+\gamma|D_1)\to1\label{interim5}
\end{equation}
for any small $\gamma>0$. Note that, conditional on $D_1$, we have $P(\xi\in\mathcal U)=P(t(\xi)\leq t(\xi^2_{(i^*)}))=F(t(\xi^2_{(i^*)}))$ where $F(\cdot)$ is the distribution function of $t(\xi)$. Since $F(t(\xi))\sim U[0,1]$ by the continuity of $t(\xi)$, we have,
\begin{eqnarray}
&&\mathbb P_{D_2}(1-\epsilon-\gamma\leq P(\xi\in\mathcal U)\leq1-\epsilon+\gamma|D_1)\label{interim8}\\
&=&P(\#\{U_i<1-\epsilon-\gamma\}\leq i^*-1,\ \#\{U_i>1-\epsilon+\gamma\}\leq n_2-i^*){}\notag\\
&&{}\text{\ \ where $\{U_i\}$ denotes $n_2$ realizations of i.i.d. $U[0,1]$ variables,}{}\notag\\
&&{}\text{\ \ $\#\{U_i<1-\epsilon-\gamma\}$ and $\#\{U_i>1-\epsilon+\gamma\}$ count the numbers of $U_i$'s that are $<1-\epsilon-\gamma$ and}\notag\\
&&{}\text{\ \ $>1-\epsilon+\gamma$ respectively}\notag\\
&\geq&1-P(\#\{U_i<1-\epsilon-\gamma\}>i^*-1)-P(\#\{U_i>1-\epsilon+\gamma\}>n_2-i^*)\label{interim7}
\end{eqnarray}
 Consider the second term in \eqref{interim7}. We have
\begin{eqnarray*}
&&P(\#\{U_i<1-\epsilon-\gamma\}>i^*-1)\\
&=&P(Bin(n_2,1-\epsilon-\gamma)>i^*-1)\\
&=&\bar\Phi\left(\frac{i^*-1-n_2(1-\epsilon-\gamma)}{\sqrt{n_2(1-\epsilon-\gamma)(\epsilon+\gamma)}}\right)+O\left(\frac{1}{\sqrt{n_2}}\right){}\\
&&{}\text{\ \ by the Berry-Essen Theorem, where $\bar\Phi$ is the tail distribution function of standard normal}\\
&=&\bar\Phi\left(\frac{i^*-1-n_2(1-\epsilon)}{\sqrt{n_2(1-\epsilon)\epsilon}}\sqrt{\frac{1-\epsilon}{1-\epsilon-\gamma}\frac{\epsilon}{\epsilon+\gamma}}+\frac{\sqrt{n_2}\gamma}{\sqrt{(1-\epsilon-\gamma)(\epsilon+\gamma)}}\right)+O\left(\frac{1}{\sqrt{n_2}}\right)\\
&\to&0\text{\ \ by \eqref{interim6}}
\end{eqnarray*}
Similarly, for the third term in \eqref{interim7}, we have
\begin{eqnarray*}
&&P(\#\{U_i>1-\epsilon+\gamma\}>n_2-i^*)\\
&=&P(Bin(n_2,\epsilon-\gamma)>n_2-i^*)\\
&=&\bar\Phi\left(\frac{n_2-i^*-n_2(\epsilon-\gamma)}{\sqrt{n_2(\epsilon-\gamma)(1-\epsilon+\gamma)}}\right)+O\left(\frac{1}{\sqrt{n_2}}\right){}\\
&&{}\text{\ \ by the Berry-Essen Theorem}\\
&=&\bar\Phi\left(-\frac{i^*-n_2(1-\epsilon)}{\sqrt{n_2(1-\epsilon)\epsilon}}\sqrt{\frac{\epsilon}{\epsilon-\gamma}\frac{1-\epsilon}{1-\epsilon+\gamma}}+\frac{\sqrt{n_2}\gamma}{\sqrt{(\epsilon-\gamma)(1-\epsilon+\gamma)}}\right)+O\left(\frac{1}{\sqrt{n_2}}\right)\\
&\to&0\text{\ \ by \eqref{interim6}}
\end{eqnarray*}
Hence \eqref{interim7} converges to 1.

\underline{Proof of \ref{thm pt1}.}
Using again the fact that, conditional on $D_1$, $F(t(\xi))\sim U[0,1]$ and $P(\xi\in\mathcal U)=F(t(\xi^2_{(i^*)}))$, we have
\begin{eqnarray*}
&&\mathbb P_{D_2}(P(\xi\in\mathcal U)\geq1-\epsilon|D_1)\\
&=&P(\#\{U_i<1-\epsilon\}\leq i^*-1)\\
&=&P(Bin(n_2,1-\epsilon)\leq i^*-1)\\
&=&\Phi\left(\frac{i^*-1-n_2(1-\epsilon)}{\sqrt{n_2(1-\epsilon)\epsilon}}\right)+O\left(\frac{1}{\sqrt{n_2}}\right)\text{\ \ by using \eqref{Berry-Essen}}\\
&\to&1-\delta\text{\ \ by \eqref{interim6}}
\end{eqnarray*}
which concludes Part \ref{thm pt1} of the theorem.
\Halmos

\endproof

Note that \eqref{interim6} is mentioned in \cite{serfling2009approximation} Section 2.6.1, and implies that, given $D_1$,
\begin{equation}
\sqrt{n_2}(P(\xi\in\mathcal U)-(1-\epsilon))=\sqrt{n_2}(F(t(\xi_{(i^*)}^2))-(1-\epsilon))\Rightarrow N\left(\sqrt{\epsilon(1-\epsilon)}\Phi^{-1}(1-\delta),\epsilon(1-\epsilon)\right)\label{new revision}
\end{equation}
by using \cite{serfling2009approximation} Corollary 2.5.2, which can be used to prove Part 1 of the theorem as well (as in \cite{serfling2009approximation} Section 2.6.3). From \eqref{new revision}, we see that $P(\xi\in\mathcal U)$ concentrates at $1-\epsilon$, as it is approximately $(1-\epsilon)+Z/\sqrt{n_2}$ where $Z\sim N\left(\sqrt{\epsilon(1-\epsilon)}\Phi^{-1}(1-\delta),\epsilon(1-\epsilon)\right)$.

\section{Illustration of Attained Theoretical Confidence Levels}\label{sec:phase 2 data}

\begin{figure}[ht]
						\centering
						\begin{minipage}[b]{0.4\textwidth}
							\centering
							\includegraphics[width=\linewidth]{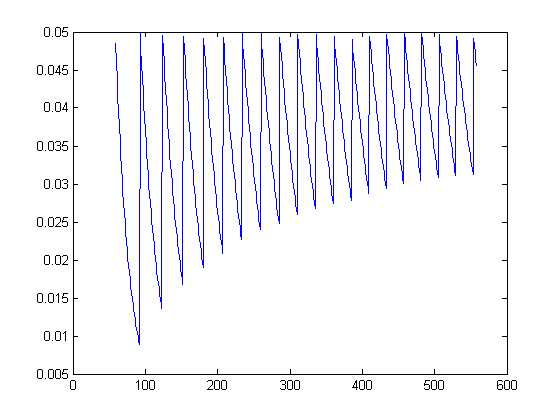}
							\caption{$\delta_{theoretical}$ against $n_2$ when $\delta=0.05$ and $\epsilon=0.05$}
							\label{fig:sub3}
						\end{minipage}%
						\begin{minipage}[b]{0.4\textwidth}
							\centering
							\includegraphics[width=\linewidth]{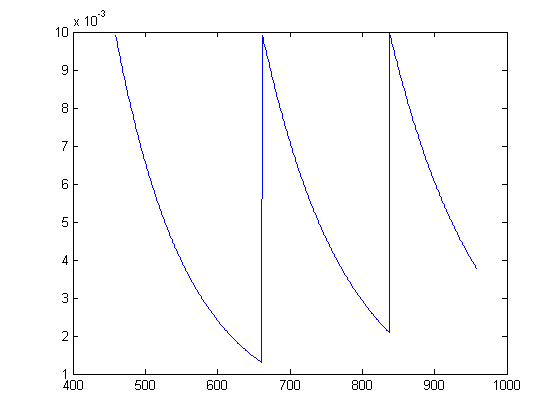}
							\caption{$\delta_{theoretical}$ against $n_2$ when $\delta=0.01$ and $\epsilon=0.01$}
							\label{fig:sub4}
						\end{minipage}
						\label{fig:f2}
				\end{figure}

The argument in Lemma \ref{quantile estimation} and the discussion after Theorem \ref{main result} implies that the theoretical confidence level for a given Phase 2 sample size $n_2$ is
		$$1-\delta_{theoretical}=\mathbb P_D(P(\xi\in\mathcal U)\geq1-\epsilon)=\sum_{k=0}^{i^*-1}\binom{n_2}{k}(1-\epsilon)^k\epsilon^{n_2-k}$$
		This quantity is in general not a monotone function of the sample size, but it does converge to $1-\delta$ as $n_2$ increases, as shown in Theorem \ref{consistency} Part \ref{thm pt1}. Figures \ref{fig:sub3} and \ref{fig:sub4} illustrate how $\delta_{theoretical}$ changes with $n_2$ for two pairs of $\epsilon$ and $\delta$. The changes follow a zig-zag pattern, with a general increasing trend. In the case $\delta=0.05$ and $\epsilon=0.05$ for example, local maxima of $\delta_{theoretical}$ occur at $n_2=59, 93, 124 ,153, 181,\ldots$
				
				\par


\section{Using RO Reformulations}\label{sec:RO results}
Results from the following discussion are adapted from \cite{bertsimas2011theory}. Further details can be found therein and in, e.g., \cite{ben2009robust}. Along with reviewing these results, we also describe how to cast them in our procedure in Section \ref{sec:procedure}.


We focus on linear safety conditions in \eqref{CCP}, i.e., $g(x;\xi)\in\mathcal A$ is in the form $Ax\leq b$, where $A\in\mathbb R^{l\times d}$ is uncertain and $b\in\mathbb R^l$ is constant. Here $A$ is identified with the random vector $\xi$. The following discussion also holds if $x$ is further constrained to lie in some deterministic set, say $\mathcal B$. For convenience, we denote each row of $A$ as $a_i'$ and each entry in $b$ as $b_i$, so that the safety condition can also be written as $a_i'x\leq b_i,i=1,\ldots,l$.


It is well-known that in solving the robust counterpart (RC), it suffices to consider uncertainty sets in the form $\mathcal U=\prod_{i=1}^l\mathcal U_i$ where $\mathcal U_i$ is the uncertainty set projected onto the portion associated with the parameters in each constraint, and so typically we consider the RC of each constraint separately.

We first consider ellipsoidal uncertainty:
\begin{theorem}[c.f. \cite{ben1999robust}]
The constraint$$a_i'x\leq b_i\ \forall a_i\in\mathcal U_i$$
where $\mathcal U_i=\{a_i=a_i^0+\Delta_iu:\|u\|_2\leq\rho_i\}$ for some fixed $a_i^0\in\mathbb R^d$, $\Delta_i\in\mathbb R^{d\times r}$, $\rho_i\in \R$, for $u\in\mathbb R^r$, is equivalent to
$${a_i^0}'x+\rho_i\| \Delta_i' x\|_2\leq b_i$$
\label{thm:ellipsoid}
\end{theorem}
											
Note that $\mathcal U_i$ in Theorem \ref{thm:ellipsoid} is equivalent to $\{a_i:\|\Delta_i^{-1}(a_i-a_i^0)\|_2\leq\rho_i\}$ if $\Delta_i$ is invertible. Thus, given an ellipsoidal set (for the uncertainty in constraint row $i$) calibrated from data in the form $\{a_i:(a_i-\mu)'\Sigma^{-1}(a_i-\mu)\leq s\}$ where $\Sigma$ is positive definite and $s>0$, we can take $a_i^0=\mu$, $\Delta_i$ as the square-root matrix in the Cholesky decomposition of $\Sigma$, and $\rho_i=\sqrt s$ in using the depicted RC.

Next we have the following result on polyhedral uncertainty:					
\begin{theorem}	[c.f. \cite{ben1999robust} and \cite{bertsimas2011theory}]
The constraint$$a_i'x\leq b_i\ \forall a_i\in\mathcal U_i$$
where $\mathcal U_i=\{a_i:D_ia_i\leq e_i\}$ for fixed $D_i\in \R^{r\times d}$, $e_i\in \R^r$ is equivalent to
	$$\begin{array}{l}
    p_i'e_i\leq b_i\\
	 		p_i'D_i=x'\\
	 		p_i\geq 0
 		\end{array}$$
where $p_i\in\mathbb R^r$ are newly introduced decision variables.\label{thm:polytope}
\end{theorem}

The following result applies to the collection of constraints $Ax \leq b$ with the uncertainty on $A\in\mathbb R^{l\times d}$ represented via a general norm on its vectorization.

\begin{theorem}[c.f. \cite{bertsimas2004robust}]
\label{thm:large_joint}
The constraint
$$Ax\leq b\ \ \forall A \in \mathcal U$$
where
	\begin{equation}
		\mathcal{U}=\{A: \|Q(vec(A)-vec(\bar{A}))\|\leq \rho\},\label{uncertainty set multiple}
	\end{equation}
for fixed $\bar A\in\mathbb R^{l\times d}$, $Q\in\R^{ld\times ld}$ invertible, $\rho\in\mathbb R$, $vec(A)$ as the concatenation of all the rows of $A$, $\|\cdot\|$ any norm, is equivalent to
$$\bar{a}'_ix+\rho \|(Q')^{-1}x_i \|^* \leq b_i,i=1,...,l$$
where $\bar{a}_i'\in\mathbb R^d$ is the $i$-th row of $\bar A$, $x_i \in \R^{(ld)\times 1}$ contains $x \in \R^d$ in entries $(i-1)d+1$ through $i\,d$ and 0 elsewhere, and $\|\cdot\|^*$ is the dual norm of $\|\cdot\|$.\label{thm:joint ellipsoid}
\end{theorem}

When $\|\cdot\|$ denotes the $L_2$-norm, Theorem \ref{thm:joint ellipsoid} can be applied in much the same way as Theorem \ref{thm:ellipsoid}, with $vec(\bar A)$ denoting the center, $Q$ taken as the square root of the Cholesky decomposition of $\Sigma^{-1}$ where $\Sigma$ is the covariance matrix, and $\rho=\sqrt s$ where $s$ is the squared radius in an ellipsoidal set constructed for the data of $vec(A)$.

Next we have the following theorem to handle \eqref{uncertainty PCA}, which can be proved similarly as for Theorem \ref{thm:ellipsoid} or by standard conic duality.

\begin{theorem}
The constraint
$$\xi'x\leq b\ \forall \xi\in\mathcal U$$
where $\mathcal U$ is defined in \eqref{uncertainty PCA}, and $\Sigma$ has full rank, is equivalent to
$$
\begin{array}{l}
\mu' \Sigma^{-1/2} u + \sqrt{s} \lambda \leq b\\
M'\Sigma^{-1/2}u=x\\
\|u\|_2\leq \lambda,
\end{array}
$$where $\lambda \in \R$, $u\in \R^ {r}$ are additional decision variables.\label{thm:PCA}
\end{theorem}

\section{Further Discussion on Choices of Uncertainty Sets}\label{sec:comparison}
This section extends the discussions in Section \ref{sec:learn tractable} on choosing suitable uncertainty sets.

\subsection{Comparing Individualized Ellipsoids and a Single Ellipsoids for Joint Chance Constraints}\label{sec:comparision details}
We state the following result that compares, in the case of joint chance constraints, between the use of individual ellipsoids for the stochasticities on different constraints and a single ellipsoid for all.

\begin{proposition}
Let $\xi\in\mathbb R^m$ be a vector that can be represented as $\xi=(\xi^i)_{i=1,\ldots,k}$ with $\xi^i\in\mathbb R^{r^i}$ and $\sum_{i=1}^kr^i=m$. Let $\mathcal U_{joint}=\{\xi:\|M(\xi-\mu)\|_2^2\leq\rho_{joint}\}$ where $M$ is a block diagonal matrix
\begin{equation}
				M=  \left( \begin{array}{cccc} M^1&&&\\ &M^2&&\\&&...&\\&&&M^k    \end{array}\right)   ,
			\end{equation}
and each $M^i\in\mathbb R^{r^i\times r^i}$. Let $\mathcal U_{individual}=\prod_{i=1}^k\mathcal U^i$ where $\mathcal U^i=\{\xi^i:\|M^i(\xi^i-\mu^i)\|_2^2\leq\rho_{individual}\}$ and $(\mu^i)_{i=1,\ldots,k}$ is defined such that $\mu=(\mu^i)_{i=1,\ldots,k}$ analogously as in $(\xi^i)_{i=1,\ldots,k}$ for $\xi$. Suppose that $\mathcal U_{joint}$ and $\mathcal U_{individual}$ are calibrated using the same Phase 2 data, with the transformation maps defined as $t_{joint}(\xi)=\|M(\xi-\mu)\|_2^2$ and $t_{individual}(\xi)=\max_{i=1,\ldots,k}\|M^i(\xi^i-\mu^i)\|_2^2$ respectively.

Consider the RO
\begin{equation}
\text{minimize\ }f(x)\text{\ \ subject to\ \ }g_i(x;\xi^i)\in\mathcal A_i,i=1,\ldots,l,\ \forall \xi\in\mathcal U\label{RO special}
\end{equation}
Let $x_{joint}$ be an optimal solution obtained by setting $\mathcal U=\mathcal U_{joint}$, and $x_{individual}$ be an optimal solution obtained by setting $\mathcal U=\mathcal U_{individual}$. We have $f(x_{joint})\geq f(x_{individual})$. In other words, using $\mathcal U_{joint}$ is more conservative than using $\mathcal U_{individual}$.\label{joint individual}
\end{proposition}

\proof{Proof of Proposition \ref{joint individual}.}
The $\rho_{joint}$ calibrated using Phase 2 data is set as $t_{joint}(\xi_{(i_{joint}^*)}^2)=\|M(\xi_{(i_{joint}^*)}^2-\mu)\|_2^2$ where $i_{joint}^*$ is defined similarly as \eqref{quantile choice}. On the other hand, the $\rho_{individual}$ in the set $\mathcal U_{individual}$ (equal among all $\mathcal U^i$), is set as $t_{individual}(\xi_{(i_{individual}^*)}^2)=\max_{i=1,\ldots,k}\|M^i(\xi_{(i_{individual}^*)}^{i,2}-\mu^i)\|_2^2$ where $(\xi_{(i_{individual}^*)}^{i,2})_{i=1,\ldots,k}$ is the corresponding partition of $\xi_{(i_{individual}^*)}^2$. Using $\|M(\xi-\mu)\|_2^2=\sum_{i=1}^k\|M^i(\xi^i-\mu^i)\|_2^2$ and the fact that $\sum_{i=1}^ky_i\geq\max_{i=1,\ldots,k}y_i$ for any $y_i\geq0$, we must have $\|M(\xi-\mu)\|_2^2\geq\max_{i=1,\ldots,k}\|M^i(\xi^i-\mu^i)\|_2^2$, and so $\rho_{joint}\geq\rho_{individual}$. Note that, when projecting to each constraint, the considered RO is written as
$$\text{minimize\ }f(x)\text{\ \ subject to\ \ }g_i(x;\xi^i)\in\mathcal A_i, \forall\xi^i\in\mathcal U^i,\ i=1,\ldots,l$$
where $\mathcal U^i=\{\xi:\|M^i(\xi^i-\mu^i)\|_2^2\leq\rho_{joint}\}$ and $\{\xi:\|M^i(\xi^i-\mu^i)\|_2^2\leq\rho_{individual}\}$ for the two cases respectively. Since $\rho_{joint}\geq\rho_{individual}$, we conclude that $f(x_{joint})\geq f(x_{individual})$.\Halmos
\endproof

Proposition \ref{joint individual} is evident in that the relation $t_{joint}(\xi)\geq t_{individual}(\xi)$ leads to a larger $\mathcal U_{joint}$ and hence a smaller resulting feasible region for \eqref{RO special} compared with $\mathcal U_{individual}$. It hints that, if the data across the constraints are uncorrelated, it is always better to use constraint-wise individual ellipsoids that are calibrated jointly. The same holds if we choose to use diagonalized ellipsoids in our representation, as these satisfy the block-diagonal structural assumption in the proposition. On the other hand, if the data across individual constraints are dependent and we want to capture their correlations in our ellipsoidal construction, the comparison between the two approaches is less clear. 

\subsection{Complexity of Uncertainty Sets}\label{sec:complexity}
Another consideration in choosing uncertainty set in our framework is the set complexity. For example, we can use an ellipsoidal set with a full covariance matrix, a diagonalized matrix and an identity matrix, the latest leading to a ball. The numbers of parameters in these sets are in decreasing order, making the sets less and less ``complex". Generally, more data supports the use of higher complexity representation, because they are less susceptible to over-fitting. In terms of the average optimal value obtained by the resulting RO, we observe the following general phenomena:
\begin{enumerate}
\item Ellipsoidal sets with full covariance matrices are generally better than diagonalized elliposids and balls when the Phase 1 data size is larger than the dimension of the stochasticity. However, if the data size is close to or less than the dimension, the estimated full covariance matrix may become singular, causing numerical instabilities. 
\item In the case where ellipsoidal sets are problematic (due to the issue above), diagonalized ellipsoids are preferable to balls unless the data size is much smaller than the stochasticity dimension.
\end{enumerate}

Note that the above observations are consistent with theoretical results in covariance matrix estimation. In particular, it is known that the data size required to accurately estimate the covariance matrix of an $m$-dimensional random vector is of order (arbitrarily higher than) $m$ if the vector is sub-Gaussian (Theorem 4.7.1 in \cite{vershynin2018high}) and $m\log m$ for more and very general vectors (Theorem 5.6.1 in \cite{vershynin2018high}). This suggests that using fully estimated covariance matrix  is desirable over diagonalized matrix when data size is slightly above the dimension.



\subsection{Missing Details for Section \ref{sec:test_case_3}}\label{sec:missing test case}
The example in Section \ref{sec:test_case_3} utilizes the observations discussed in Appendices \ref{sec:comparision details} and \ref{sec:complexity}, which we detail below. 
%
Since the sample size is less than the stochasticity dimension, we use diagonalized ellipsoids in our constructions. Next, we compare using individualized ellipsoids each for the stochasticity in each constraint versus a single ellipsoid, as depicted in Proposition \ref{joint individual}. Table \ref{table:joint_ccp_compare} column 2 shows the results using a single ellipsoid over vectorized $A$, and column 3 shows the counterparts for individually constructed ellipsoids. We observe that the latter has a smaller average optimal value (-6957.26 versus -4529.51), which is consistent with the implication from Proposition \ref{joint individual}. 

\begin{table}[t]
	\centering
	\caption{Comparing the optimality and feasibility performances between single diagonalized ellipsoid and individually constructed diagonalized ellipsoids, under sample size $n=120$, and we use $n_1=60$ and $n_2=60$.}
	\label{table:joint_ccp_compare}
	\begin{tabular}{|l|l|l|}
		\hline
		& RO(Single Diagonalized Ellipsoid)       & RO(Individual Diagonalized Ellipsoids)     \\ \hline
		Obj. Val.      & -4529.51 & -6957.26  \\ \hline
		$\hat{\epsilon}$ & 0 & $3.55\times 10^{-5}$            \\ \hline
		$\hat{\delta}$     & 0    & 0             \\ \hline
	\end{tabular}
\end{table}

\begin{table}[t]
	\centering
	\caption{Comparing the optimality and feasibility performances between two scaling strategies for reconstructing the uncertainty set.}
	\label{table:joint_ccp_scale}
	\begin{tabular}{|l|l|l|}
		\hline
		& Reconstructed RO (Scale 1)       & Reconstructed RO (Scale 2)      \\ \hline
		Obj. Val.      & -7880.06 & -7541.29  \\ \hline
		$\hat{\epsilon}$ & 0.0127 & 0.0017           \\ \hline
		$\hat{\delta}$     & 0.029    & 0             \\ \hline
	\end{tabular}
\end{table}

We further investigate the use of reconstruction for joint CCP. We use $\max_{j=1,...,l} \{(a_j'\hat{x}_0-b_j)/k_j \}$ to determine the quantile for calibrating the size of the uncertainty set, where $k_j$ is a scale parameter assigned to constraint $j$. Table \ref{table:joint_ccp_scale} compares two natural choices of $k_j$ for the same problem as above but with a different $\Sigma$. Column 2 uses $k_j=b_j-\mu_j'\hat{x}_0$, where $\mu_j'$ is the sample mean of the Phase 1 data of $a_j'$. Column 3 uses $k_j=std(a_j'\hat{x}_0)$, the standard deviation of the Phase 1 data of $a_j'\hat{x}_0$. While the performances using these two scale parameters can be problem dependent, we observe that the former works better in this example (with a better average optimal value) and hence adopt it for our experiment.

\section{Integrating with Machine Learning Tools} \label{sec:single_linear_ml}
We provide some numerical results to support the use of the machine learning tools described in Section \ref{sec:learn tractable}. Throughout this section we use the single CCP (\ref{eq:single_CCP_ex}) as an example.

\subsection{Cluster Analysis}

To illustrate the use of clustering, suppose $\xi$ follows a mixture of $N(\mu_1, \Sigma_1)$ and $N(\mu_2, \Sigma_2)$ with probabilities $\pi_1=\pi_2=0.5$. Table \ref{cluster} column 2 shows the performance of our RO using a single ellipsoidal set. Column 3 shows the result when we first apply 2-mean clustering to Phase 1 data and construct a union of ellipsoids. The average objective value (-961.434) is demonstrably improved compared to using a single ellipsoid (-940.502). Similarly, the reconstructed RO from using clustering performs better than RO using a single ellipsoid, and both are better than the non-reconstructed counterparts. 


\begin{table}[t]
	\centering
	\caption{Optimality and feasibility performances on a single linear CCP with mixture Gaussian distributions for several methods, under sample size 300, and we use $n_1=240$ and $n_2=60$.}
	
	\resizebox{\columnwidth}{!}{%
	\begin{tabular}{|l|l|l|l|l|}
		\hline
		&  RO(Unclustered) & RO(Clustered) & Reconstructed RO(Unclustered) & Reconstructed RO(Clustered)       \\ \hline
		Obj. Val.    &  -940.502    & -961.434 & -1074.63 & -1087.66  \\ \hline
	    $\hat{\epsilon}$ & $2.18\times 10^{-7}$    & $3.01\times 10^{-6}$ & 0.0162 & 0.0163   \\ \hline
		
		$\hat{\delta}$  & 0           & 0        & 0.05     & 0.049     \\ \hline
		
	\end{tabular}}

	\label{cluster}
\end{table}

\subsection{Dimension Reduction}\label{sec:dimension reduction}

To illustrate the use of dimension reduction, we specify $\xi$ as follows. We first generate  $\tilde{\xi} \in \R^{11}$ under $N(\mu, \Sigma)$, where $\mu$ and $\Sigma$ are arbitrary vector and positive definite matrix. We create a higher dimensional $\xi\in \R^{1100}$ by $\xi=P\tilde{\xi}+\omega$, where $\omega$ is a ``perturbation" vector with each element distributed uniformly on [-0.0005,0.0005] and $P \in \R ^{1100 \times 11}$.

Table \ref{table:pca_high_dim} column 2 shows the results using RO with a diagonalized ellipsoid on the data of $\xi$. Diagonalized ellipsoid is used here because the dimension $d=1100$, which is much larger than the Phase 1 data size $n_1=60$, causes singularity issue when constructing a full ellipsoid. Column 3 shows the results when we apply principal component analysis (PCA) to reduce the data to the 11 components having the largest variances and use the linearly transformed ellipsoid \eqref{uncertainty PCA}. The number of components 11 is chosen from the cutoff of leaving out $0.01\%$ of the total variance, which we declare as negligible. The PCA approach outperforms the use of a basic diagonalized ellipsoid in terms of average optimal value (-1189.32 versus -1039). 


As can be seen in this example, the dimension reduction brought by PCA allows to use a full ellipsoid that captures the shape of the data better on the relevant directions than using the original data, whose high dimension forces one to adopt a simpler geometric set such as diagonalized ellipsoid. Our recommendation in selecting the number of components in PCA is to be conservative, in the sense of choosing one as large as possible so long as it is small enough to support the use of a full ellipsoid (roughly speaking, this means it is smaller than the Phase 1 data size). 

\begin{table}[t]
	\centering
	\caption{Optimality and feasibility performances on a $d=1100$ dimensional single linear CCP using PCA, under sample size $n=120$, and we use $n_1=60$ and $n_2=60$.}
	\label{table:pca_high_dim}
	\begin{tabular}{|l|l|l|}
		\hline
		&  RO(Diagonalized Ellipsoid) & RO(PCA with 11 Components)     \\ \hline
		Obj. Val.       & -1039    & -1189.32  \\ \hline
		$\hat{\epsilon}$    & $4.54 \times 10^{-16}$        & $1.43 \times 10^{-5}$        \\ \hline
		$\hat{\delta}$        & 0      & 0        \\ \hline
	\end{tabular}
\end{table}

\subsection{``Basis" Learning}\label{sec:basis}

We consider the last approach described in Section \ref{sec:learn tractable} that surrounds each observation with a ball. For convenience, we call this approach ``basis" learning (as we view each of these created balls as a ``basis"). We set $\xi \sim \mathcal{N}(\mu, \Sigma)$ for some arbitrarily chosen $\mu$ and $\Sigma$ and $d=11$. Table \ref{table:basis_set} shows that the basis learning approach (column 4) outperforms the use of a diagonalized ellipsoid (column 3), but underperforms the use of a full ellipsoid (column 2), in terms of average optimal value (-1016.95, -946.33 and -1186.86 respectively). All three approaches are conservative however ($\hat{\delta}\approx0$). This roughly indicates that basis learning is capable of capturing some covariance information. 
\begin{table}[t]
	\centering
	\caption{Optimality and feasibility performances on a single linear CCP for basis learning and other methods, under sample size $n=80$, and we use $n_1=21$ and $n_2=59$.}
	\label{table:basis_set}
	\begin{tabular}{|l|l|l|l|}
		\hline
		&  RO(Ellipsoid)& RO(Diagonalized Ellipsoid) & RO(Basis)     \\ \hline
		Obj. Val.      & -1186.86  &-946.33  & -1016.95  \\ \hline
		$\hat{\epsilon}$    & 0.0002  &$3.03 \times 10^{-4}$      & $1.22 \times 10^{-8}$        \\ \hline
		$\hat{\delta}$        & 0      & 0  &0      \\ \hline
	\end{tabular}
\end{table}

Next we generate $\xi$ from a mixture of Gaussian distribution with 5 components and $d=11$. Table \ref{table:basis_cluster} shows that basis learning (column 4) outperforms ellipsoid (column 2) in terms of average optimal value (-1033.84 versus -845.973). However, it does not perform as well compared to using the union of 5 ellipsoids from clustering (column 3, with an average optimal value -1090.57). This supports the guidance that,  when applying to convoluted data, basis learning is better than using over-simplified shape, but may not work as well compared to other established machine learning tools.


\begin{table}[t]
	\centering
	\caption{Optimality and feasibility performances on a single linear CCP for basis learning and other methods, using sample size $n=300$. For learning-based RO, we use $n_1=240$ and $n_2=60$.}
	\label{table:basis_cluster}
	\begin{tabular}{|l|l|l|l|l|}
		\hline
        
		& SG       & RO(Ellipsoid) & RO(Clustered) & RO(Basis)    \\ \hline
		Obj. Val.      & -1191.82 & -845.973    & -1090.57 & -1033.84  \\ \hline
		$\hat{\epsilon}$  & 0.037   & $2.20 \times 10^{-5}$        & $8.73 \times 10^{-12}$    & 0      \\ \hline
		$\hat{\delta}$     & 0.125   & 0      & 0   & 0       \\ \hline
	\end{tabular}
\end{table}

\section{Tractable Reformulation of DRO under Ellipsoidal Moment-Based Uncertainty Set}\label{sec:quadratic DRO}
We review the tractable reformulation of moment-based DRO. In particular, we focus on the extension of the DRO reformulation under first and second moment information in \cite{ghaoui2003worst} using the ellipsoidal uncertainty set suggested in \cite{marandi2017extending}. 

For single linear CCP with constraint $P(\xi'x\leq b)\geq 1-\epsilon$,  \cite{ghaoui2003worst} shows that the worst-case constraint, among all distributions generating $\xi$ that have exactly known mean $\mu$ and covariance matrix $\Sigma$, can be reformulated as \begin{equation}\label{eq:dro_true}
    \sqrt{\frac{1-\epsilon}{\epsilon}} \|\Sigma^{\frac{1}{2}}x\|_2-\mu' x-b\leq 0.
\end{equation}
In the situation where $\mu$ and $\Sigma$ are unknown but i.i.d. data are available, we can construct an ellipsoidal moment set $\mathcal V$ such that $P\left( (\mu,\Sigma) \in \mathcal V \right) \geq 1- \delta$, using the delta method in Section 5 of \cite{marandi2017extending}. We then consider the worst-case chance constraint over distributions with mean and covariance matrix inside $\mathcal V$, i.e., 
\begin{equation}
\inf_{Q:(E_Q[\xi],E_Q[(\xi-E_Q[\xi])(\xi-E_Q[\xi])'])\in\mathcal V} Q(\xi'x\leq b)\geq1-\epsilon\label{DRO moment formulation appendix}
\end{equation}
where $Q$ is a distribution generating $\xi\in\mathbb R^d$, $E_Q[\xi]$ is the mean and $E_Q[(\xi-E_Q[\xi])(\xi-E_Q[\xi])'])$ the covariance matrix under $Q$. Given \eqref{eq:dro_true}, the following theorem that extends the result in \cite{marandi2017extending} can be used to provide a tractable reformulation for this worst-case chance constraint.

\begin{theorem}
Let  $u \in \R$, $\hat \Gamma \in \R^{d\times d}$, $\hat w \in \R^d$, $B \in \R^{\frac{d^2+3d}{2}\times \frac{d^2+3d}{2}}$, $\rho \in \R$ be given. We set $svec(\Gamma)=[\Gamma_{11},\sqrt{2}\Gamma_{12},...,\sqrt{2}\Gamma_{1n},\Gamma_{22},...,\sqrt{2}\Gamma_{(n-1)n},\Gamma_{nn}]'$.
The constraint
\begin{equation}\label{eq:dro_ro}
    \sqrt{x'\Gamma x}+w'x+u\leq0, \forall \left( \begin{array}{c}
    w\\svec(\Gamma)
\end{array}\right)\in \mathcal U
\end{equation}
with decision variable $x\in \R^d$, where $\mathcal U=\mathcal U_1 \cap \mathcal U_2$ and $$\mathcal U_1=\left\{ \left( \begin{array}{c}
    w\\svec(\Gamma)
\end{array}\right) = B\nu+\left(\begin{array}{c}
    \hat w\\svec(\hat \Gamma)
\end{array}\right): \|\nu\|_2 \leq \rho, \nu \in \R ^{\frac{d^2+3d}{2}}\right\},
$$
$$\mathcal U_2=\left\{ \left( \begin{array}{c}
    w\\svec(\Gamma)
\end{array}\right): w\in \R^d,\ \Gamma \in S^{+}_d\right\},$$ is equivalent to  \begin{equation} \label{eq:data_dro}
    \hat w'x+ trace(\hat \Gamma W)+ \rho\left\Vert B'\left( \begin{array}{c}
    x\\svec(W)
\end{array}\right) \right\Vert_2+u+\frac{\eta}{4}\leq 0, \left[ \begin{array}{cc}
    W&x\\x'&\eta 
\end{array}\right] \succeq 0_{(d+1)\times (d+1)}
\end{equation}
where $W \in \R^{d\times d}$ and $\eta\in \R$ are additional (dummy) variables, and $ 0_{(d+1)\times (d+1)}$ is a zero matrix of size $(d+1)\times(d+1)$.

\label{thm:dro_moment_ci}
\end{theorem}

Theorem \ref{thm:dro_moment_ci} is an application of Theorem 1 (II) in \cite{marandi2017extending} on ellipsoidal uncertainty sets in the form of $\mathcal U$. 
Note that $\mathcal U$ consists of two intersecting sets, the ellipsoidal set $\mathcal U_1$ constructed from the delta method discussed in \cite{marandi2017extending} that is designed to contain the true moments of $\xi$ with confidence $1-\delta$, and the set $\mathcal U_2$ that constrains the covariance matrix to be positive semidefinite. We reformulate the worst-case chance constraint \eqref{DRO moment formulation appendix} into a semidefinite constraint by rewriting the former in the form \eqref{eq:dro_ro} using \eqref{eq:dro_true} and applying Theorem \ref{thm:dro_moment_ci}. 

When $\xi$ has dimension $d$, the total number of the first and second moments is $(3d+d^2)/2$. To form an ellipsoidal set for all these moments using the delta method, one would need to use the estimated covariance matrix for all these moments, which requires estimating higher-order moments and has size $(3d+d^2)/2 \times (3d+d^2)/2$ (for more details, see Section 5 of \cite{marandi2017extending}). The resulting optimization problem is a semidefinite program with $(5d+3d^2)/2+1$ decision variables. 

\section{Additional Numerical Results}
\label{sec:numerics added}

This section shows three additional sets of numerical results. The first is the same example as Section \ref{sec:test_case_1} but with additional non-negativity constraints. These constraints are added to make sure that SG and FAST do not generate unbounded solutions. The second set of results contain a random right hand side quantity in a linear chance constraint. It illustrates how one can use our reconstruction to enhance performance by transforming the safety condition, in the case that a direct use seems un-usable at first. Lastly, we present some further numerical investigation of divergence-based DRO. 

\subsection{Multivariate Gaussian on a Single Chance Constraint with Non-negativity Conditions}
\label{sec:append_single_example}
We consider a modification of the example in Section \ref{sec:test_case_1}
\begin{equation}
	\text{minimize\ }c'x\text{\ \ subject to\ \ }P(\xi'x\leq b)\geq1-\epsilon,\ x\geq 0
	\label{eq:single_CCP_bounded}
\end{equation}
where we add a non-negativity constraint and keep all other parts unchanged. We again consider $d=11$ and $100$. The main purpose of the modification is to eliminate the unbounded solutions that occurred in the $d=100$ case of (\ref{eq:single_CCP_ex}) when we apply SG and FAST. The comparisons among different approaches on this problem, shown in Tables \ref{table.n=11.bounded.less}, \ref{table.n=11.bounded.more}, \ref{table.n=100.bounded.less} and \ref{table.n=100.bounded.more}, are largely similar to those in Section \ref{sec:test_case_1}, but also bear some notable differences that we highlight here.

\begin{table}
	\centering
	\caption{Optimality and feasibility performances on a single $d=11$ dimensional linear CCP with non-negativity constraints for several methods, using sample size $n=120$. The true optimal value is -1106.23.}

        \begin{tabular}{|l||l|l|l|l|l|l||l|l|l|l||l|l|}
\hline
                   & RO                   & Recon     & SG         & FAST         & DRO Mo               & DRO KL    & SCA      \\ \hline
$n$                & 120                   & 120      & 120         & 120          & 120                   & 120       & -        \\ \hline
$n_1$              & 60                    & 60       & -           & 61           & -                      & 60        & -        \\ \hline
$n_2$              & 60                    & 60       & -           & 59           & -                     & 60       & -        \\ \hline
Obj. Val.          & -924.05              & -1070.75  & -1068.17  & -1060.04      & -893.84              & 0 & -1065.59 \\ \hline
$\hat{\epsilon}$   &  $4.99\times 10^{-7}$ & 0.0158   & 0.0155    & 0.0119      & $8.46\times 10^{-10}$ & 0 & 0.0072   \\ \hline
$\hat{\delta}$     & 0                     & 0.032    & 0.019       & 0.008        & 0                      & 0        & 0       \\ \hline
\end{tabular}
	
	\label{table.n=11.bounded.less}
\end{table}

\begin{table}
	\centering
	\caption{Optimality and feasibility performances on a single $d=11$ dimensional linear CCP with non-negativity constraints for several methods, using sample size $n=336$. The true optimal value is -1106.23.}
	\resizebox{\columnwidth}{!}{%

        \begin{tabular}{|l||l|l|l|l|l|l||l|l|l|l||l|l|}
\hline
                 & RO                    & Recon     & SG       & FAST      & DRO Mo             & DRO KL     & SCA      \\ \hline
$n$               & 336                    & 336     & 336        & 336       & 336                 & 336       & -        \\ \hline
$n_1$            & 212                    & 212      & -         & 318       & -                   & 168        & -        \\ \hline
$n_2$             & 124                    & 124     & -          & 18        & -                   & 168       & -        \\ \hline
Obj. Val.         & -956.63               & -1086.28 & -1050.52 & -1049.82   & -921.232            & 0 & -1065.59 \\ \hline
$\hat{\epsilon}$   & $1.34\times 10^{-6}$ & 0.0244   &  0.00534 & 0.00523 & $6.15\times 10^{-9}$&  0 & 0.0072   \\ \hline
$\hat{\delta}$    & 0                      & 0.045   & 0          & 0         & 0                   & 0       & 0        \\ \hline
\end{tabular}
	}
	
	\label{table.n=11.bounded.more}
\end{table}

\begin{table}
	\centering
	\caption{Optimality and feasibility performances on a single $d=100$ dimensional linear CCP with non-negativity constraints for several methods, using sample size $n=120$. The true optimal value is -1195.3. Results on moment-based DRO are based on 30 replications due to high computational demand.}

        \begin{tabular}{|l||l|l|l|l|l|l||l|l|l|l||l|l|}
\hline
                 & RO      & Recon    & SG       & FAST   & DRO Mo       & DRO KL     & SCA      \\ \hline
$n$              & 120      & 120     & 120      & 120     & 120           & 120       & -      \\ \hline
$n_1$            & 60       & 60      & -        & 61      & -          & 60        & -      \\ \hline
$n_2$            & 60       & 60      & -        & 59      & -         & 60       & -      \\ \hline
Obj. Val.        & -832.142 & -1111.04& -1195.26 & -980.64  & -1120.37   &  0 & -1152.35 \\ \hline
$\hat{\epsilon}$ & 0        & 0.0159  & 0.458    & 0.0170  & 0.0095     &  0  & 0.0072   \\ \hline
$\hat{\delta}$   & 0        & 0.046   & 1        & 0.064   & 0           & 0        & 0        \\ \hline
\end{tabular}
	
	\label{table.n=100.bounded.less}
\end{table}

\begin{table}
	\centering
	\caption{Optimality and feasibility performances on a single $d=100$ dimensional linear CCP with non-negativity constraints for several methods, using sample size $n=2331$. The true optimal value is -1195.3. Results on moment-based DRO are based on 30 replications due to high computational demand.}

        \begin{tabular}{|l||l|l|l|l|l|l||l|l|l|l||l|l|}
\hline
                 & RO        & Recon         & SG               & FAST     & DRO Mo               & DRO KL     & SCA      \\ \hline
$n$               & 2331        & 2331       & 2331             & 2331     & 2331                 & 2331        & -      \\ \hline
$n_1$            & 1318        & 1318        & -                & 2326     & -                    & 1166        & -      \\ \hline
$n_2$            & 1013        & 1013        & -                & 5        & -                    & 1165       & -      \\ \hline
Obj. Val.        & -1005.62   & -1164.47       & -1156.76       & -1155.51 & -1033.58             & 0& -1152.35 \\ \hline
$\hat{\epsilon}$ & 0          & 0.0397        & 0.0293          & 0.0272   & $5.18\times 10^{-11}$ & 0  & 0.0072 \\ \hline
$\hat{\delta}$    & 0          & 0.058        & 0               & 0        & 0                     & 0        & 0       \\ \hline
\end{tabular}
	
	\label{table.n=100.bounded.more}
\end{table}

In the $d=100$ case, when sample size is small ($n=120$), SG and FAST can now obtain bounded solutions. However, SG fails to obtain feasible solutions as shown by $\hat\delta=1$ in Table \ref{table.n=100.bounded.less}, because the sample size is far smaller than the minimum requirement (2331). FAST obtains confidently feasible solutions that perform better in objective value than our plain RO (-980.64 versus -832.142), but worse than our reconstructed RO (-1111.04), the latter plausibly attributed to the initial solutions of FAST that are not in good quality.


In the $d=11$ case, SG now achieves feasibility with $n=120$ samples, and when the minimum required sample size $n=336$ is used, the solution appears more conservative compared to the counterpart in Section \ref{sec:test_case_1}, as shown by $\hat\delta=0$ in Table \ref{table.n=11.bounded.more} versus $\hat\delta=0.056$ in Table \ref{table.n=11.more}. This can be explained by the obtained solutions in the current problem being non-fully-supported (i.e., the number of support constraints is less than $d$, which gives the problem a lower ``intrinsic'' dimension). Note that when the sample size increases from 120 to 336, the solutions of SG necessarily become more conservative (regardless of the dimension in consideration), which is a consequence of the nature of constraint addition in SG.
On the other hand, the solutions in our RO improve as sample size increases, plausibly attributed to a better estimation of HPR. Reconstructed RO provides better solutions than SG and FAST in all four sets of experiments. Nonetheless, we should mention that some constraint removal approaches like sampling-and-discarding in \cite{campi2011sampling} are available to enhance the performances of SG. Finally, since the performances of DROs and SCA follow similarly as in Section \ref{sec:test_case_1}, we do not restate the comparisons with them here.

\subsection{Multivariate Gaussian on a Single Chance Constraint with Random Right Hand Side}
\label{sec:apped_rhs_random_case}

We continue to consider the single linear CCP in \eqref{eq:single_CCP_ex}, but with the right hand side quantity $b$ being random. Specifically, we set $b$ to be generated from a Gaussian distribution with mean  $1200$ and variance $100$ (in this case, $b$ is almost positive for sure). The rest of the problem follows from Section \ref{sec:test_case_1}. Note that, by the discussion at the end of Section \ref{sec:theory}, a direct use of reconstruction would not improve the solution in this example. However, we can divide $b$ on both sides of the inequality in the safety condition, which now gives a right hand side value 1 and transformed stochastiticities as the ratios of $\xi$ and $b$.

Tables \ref{table.random.rhs.less} and \ref{table.random.rhs} present the experiments on a $d=11$ dimensional problem with $n=120$ and $n=336$ sample sizes respectively. The performances of the presented approaches are consistent with the experiments in Sections \ref{sec:test_case_1} and \ref{sec:test_case_2}. Specifically, when the sample size is small ($n=120$), our RO is preferable to SG, as it obtains feasible solutions while SG fails. Reconstruction applied on the described transformed problem continues to work and perform competitively against FAST and SCA. In particular, when $n=120$, it outperforms FAST in terms of achieved objective value, but slightly falls short of SCA. When $n=336$, reconstructed RO, SG, FAST and SCA all perform very similarly. Note that SCA have assumed moment information and hence are given an upper hand in this example.

DROs contine to be conservative in this experiment. Moment-based DRO is outperformed by both plain and reconstructed ROs in both the $n=120$ and $n=336$ cases. Similar to the example in Section \ref{sec:test_case_1}, KL-divergence-based DRO obtains an adjusted tolerance level $\epsilon^*=0$, which forces the decision $x$ to satisfy the safety condition $\xi'x \leq b$ for all $\xi\in\R^d, b\in \R$, and in this case leads to an infeasible problem. 


\begin{table}[t]
	\centering
	\caption{Optimality and feasibility performances on a single $d=11$ dimensional linear CCP with random right hand side for several methods, using sample size $n=120$.}
        
        \begin{tabular}{|l||l|l|l|l|l|l||l|l|l|l||l|l|}
\hline
                 & RO        & Recon         & SG               & FAST       & DRO Mo         & DRO KL     & SCA      \\ \hline
$n$              & 120      & 120     & 120      & 120   & 120             & 120        & -      \\ \hline
$n_1$            & 60       & 60      & -        & 61                & -      & 60     & -      \\ \hline
$n_2$            & 60       & 60      & -        & 59               & -       & 60    & -      \\ \hline
Obj. Val.        & -1143.45   & -1173.62       & -1182.90       & -1167.61     & -1138.49 & infeasible & -1175.05 \\ \hline
$\hat{\epsilon}$ & $7.60\times 10^{-6}$     & -     & 0.0170        & 0.0910           & $1.00 \times 10^{-7}$  & - & 0.0074 \\ \hline
$\hat{\delta}$    & 0          & 0.045        & 0.958               & 0.053             & 0    & -        & 0        \\ \hline
\end{tabular}
	
	\label{table.random.rhs.less}
\end{table}

\begin{table}[t]
	\centering
	\caption{Optimality and feasibility performances on a single $d=11$ dimensional linear CCP with random right hand side for several methods, using sample size $n=336$.}
        
        \begin{tabular}{|l||l|l|l|l|l|l||l|l|l|l||l|l|}
\hline
                 & RO        & Recon         & SG               & FAST  & DRO  Mo             & DRO KL    & SCA      \\ \hline
$n$               & 336                    & 336     & 336        & 336                       & -   &  336    & -        \\ \hline
$n_1$            & 212                    & 212      & -         & 318              & -     & 168   & -        \\ \hline
$n_2$             & 124                    & 124     & -          & 18             & -      & 168    & -        \\ \hline
Obj. Val.        & -1149.31   & -1175.70       & -1178.00       & -1178.01     & -1143.70 & infeasible & -1175.05 \\ \hline
$\hat{\epsilon}$ & $1.60\times 10^{-6}$          & 0.0253        & 0.051           & 0.0238  & $1.00 \times 10^{-7}$ & -   & 0.0074 \\ \hline
$\hat{\delta}$    & 0          & 0.035        & 0.051               & 0.052          & 0       & -        & 0       \\ \hline
\end{tabular}
	
	\label{table.random.rhs}
\end{table}

\subsection{Additional Numerical Investigation on DRO with KL Divergence}\label{sec:phi_investigation}

We provide more details on constructing KL-divergence balls in DRO, which has been used in our numerical comparisons. In the case of continuous distributions for generating $\xi$, constructing KL balls requires estimating a reference distribution $f_0$ (center of the ball) using kernel density estimation, and then a $k$-NN or other similar methods to estimate the set size. This selection of the reference distribution aims to approximate the true distribution as much as possible, and the set size is chosen such that the divergence ball contains the true distribution with high confidence. Below we detail these procedures, followed by a very low-dimensional example where these procedures work in calibrating DRO and allow illustrative comparisons with other approaches.

\subsubsection{Bandwidth Selection for Kernel Density Estimation}

Following \cite{jiang2016data}, we use kernel density estimation to estimate the reference distribution $f_0$. This estimation procedure requires the proper selection of a bandwidth parameter, whose theoretical optimal choice is of order $N^{-\frac{1}{m+4}}$, where $N$ is the sample size and $m$ is the dimension of the randomness. In the following, we consider bandwidth in the form of $BN^{-\frac{1}{m+4}}$ for some $B\in \mathbb{R}$. 
 
\begin{figure}[t]
							\centering
							
							\begin{minipage}[b]{0.4\textwidth}
								\centering
								\includegraphics[width=\linewidth]{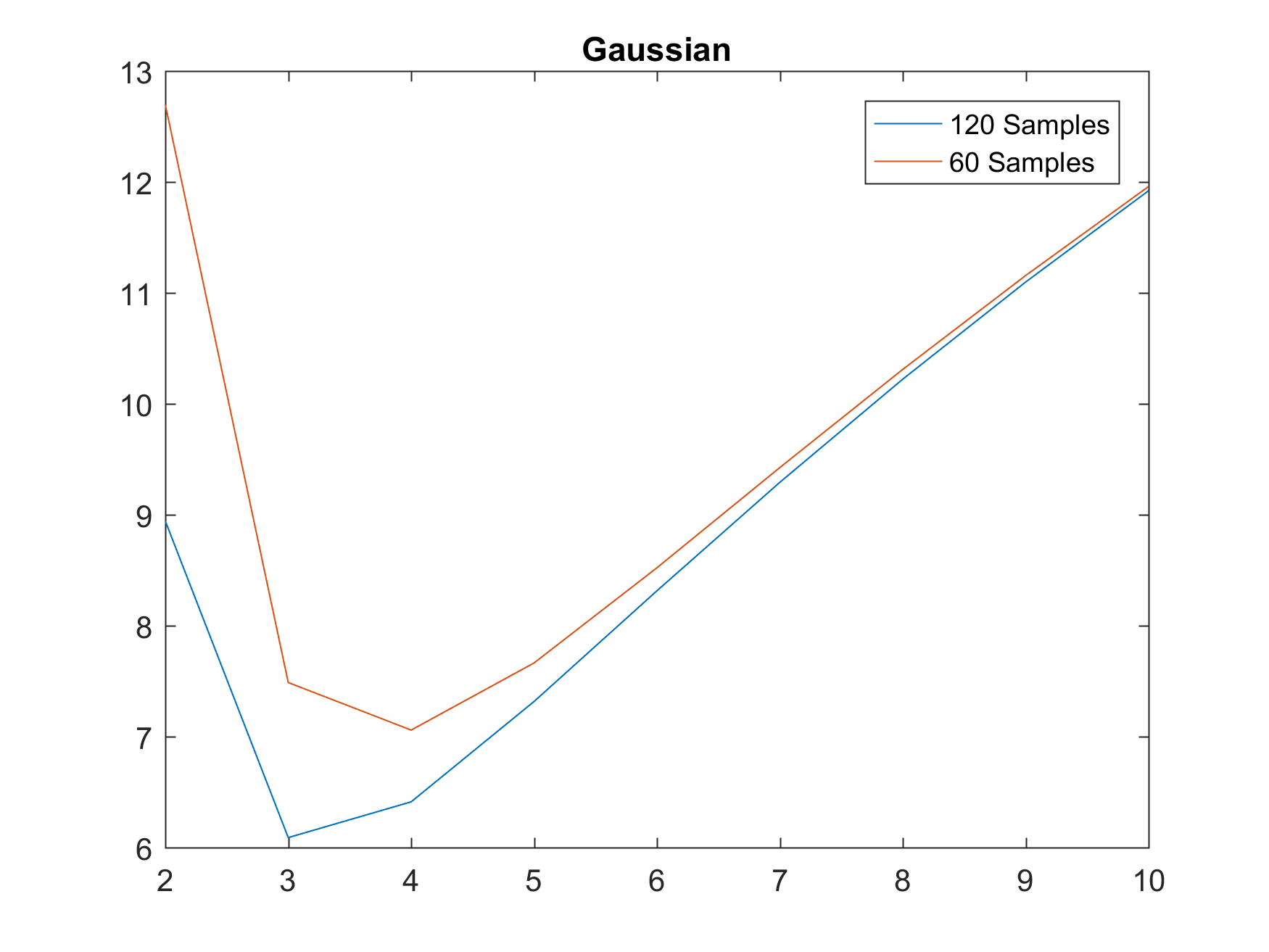}
								\caption{Divergence with different bandwidth parameter $B$ and sample size $N=120,60$. The randomness is Gaussian distributed with dimension $m=11$.}
								\label{fig:bandwitdh120}
							\end{minipage}
                            \begin{minipage}[b]{0.4\textwidth}
								\centering
                            \includegraphics[width=\linewidth]{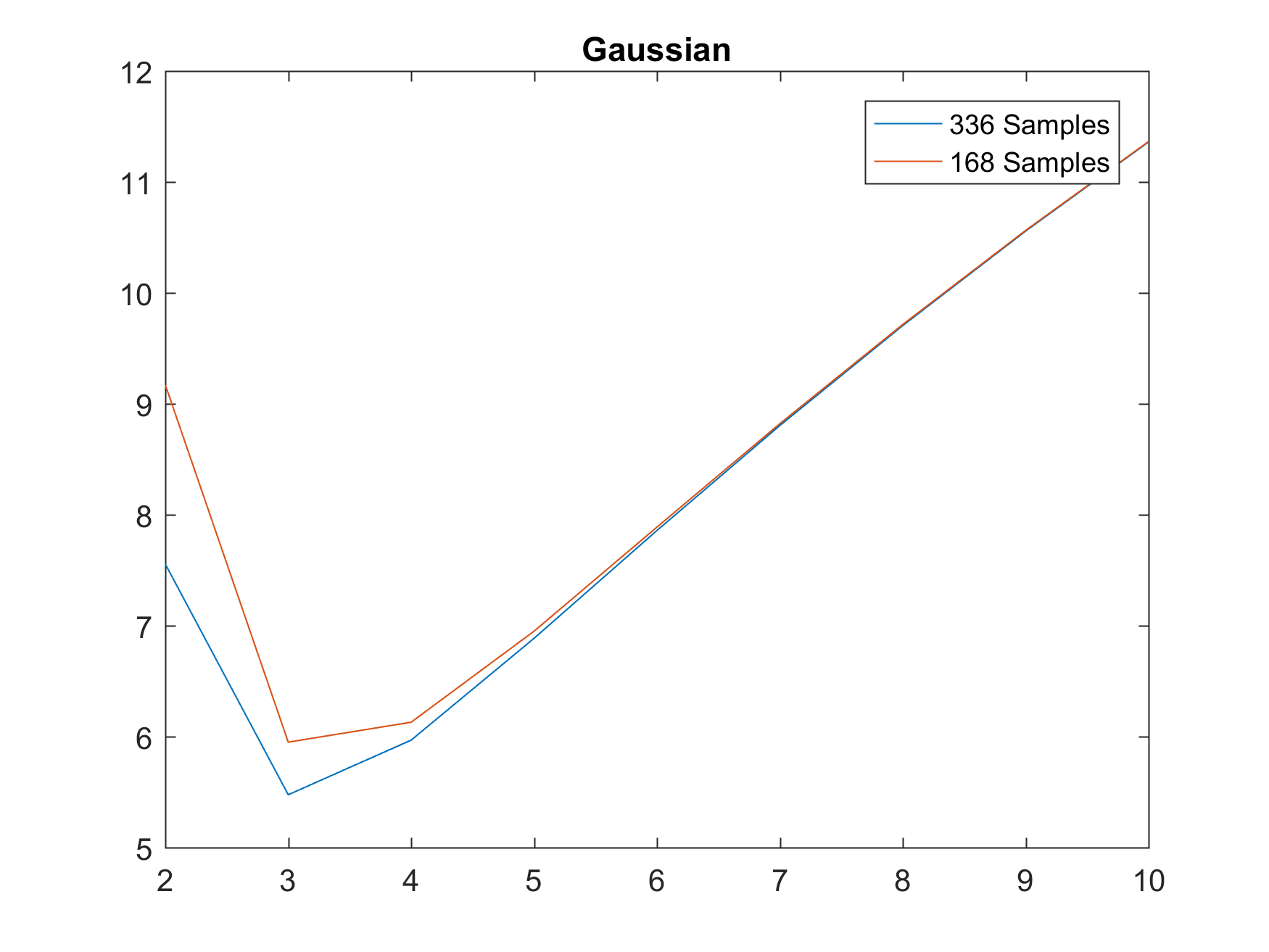}
	\caption{Divergence with different bandwidth parameter $B$ and sample size $N=336,168$. The randomness is Gaussian distributed with dimension $m=11$.}
	\label{fig:bandwitdh336}
							\end{minipage}\vspace*{-2ex}%
			\end{figure}

\begin{figure}[t]
							\centering
							
							\begin{minipage}[b]{0.4\textwidth}
								\centering
								\includegraphics[width=\linewidth]{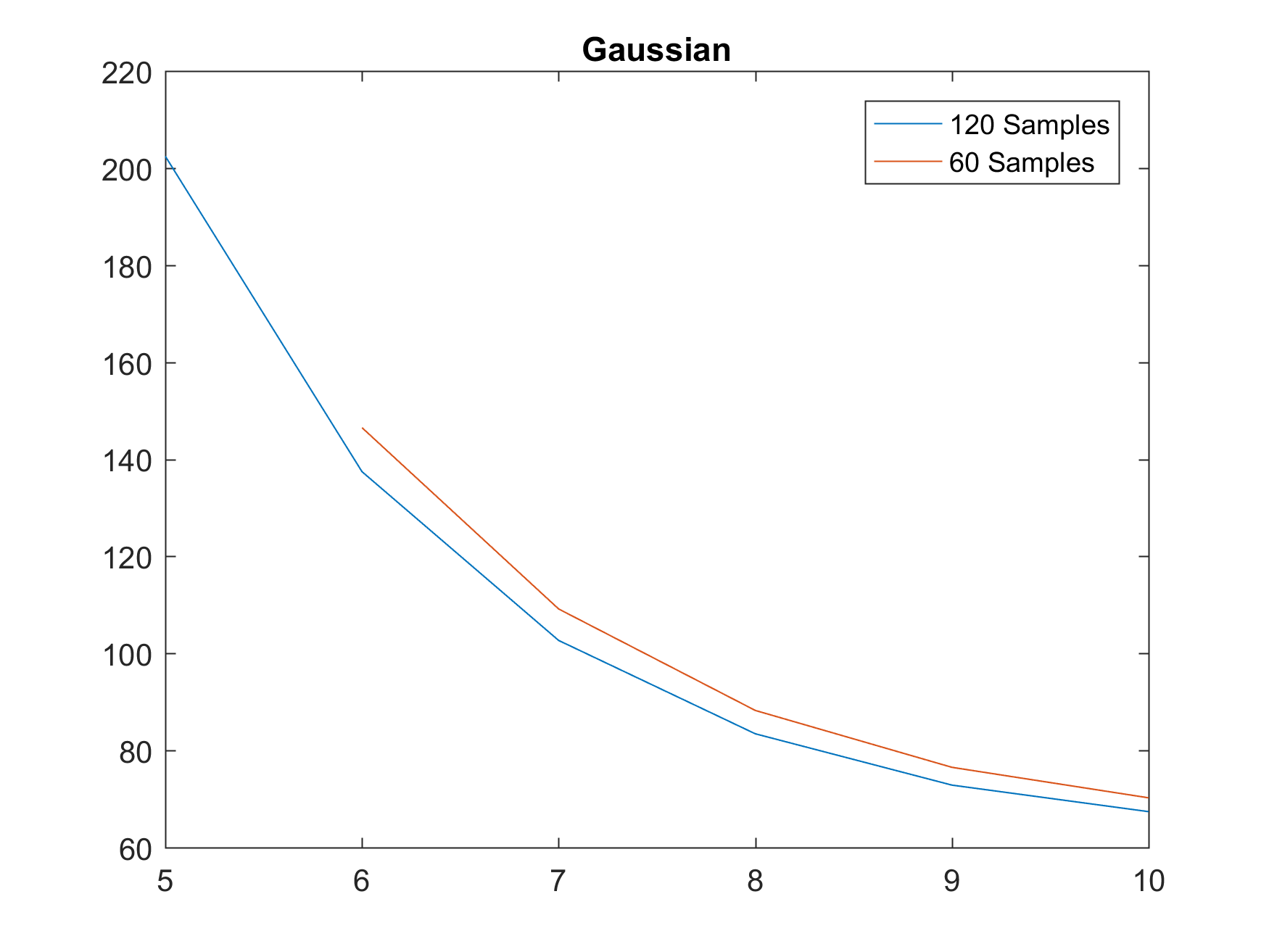}
								\caption{Divergence with different bandwidth parameter $B$. The randomness is Gaussian distributed with dimension $m=100$.}
								\label{fig:bandwitdh120d10}
							\end{minipage}
                            \begin{minipage}[b]{0.4\textwidth}
								\centering
                            \includegraphics[width=\linewidth]{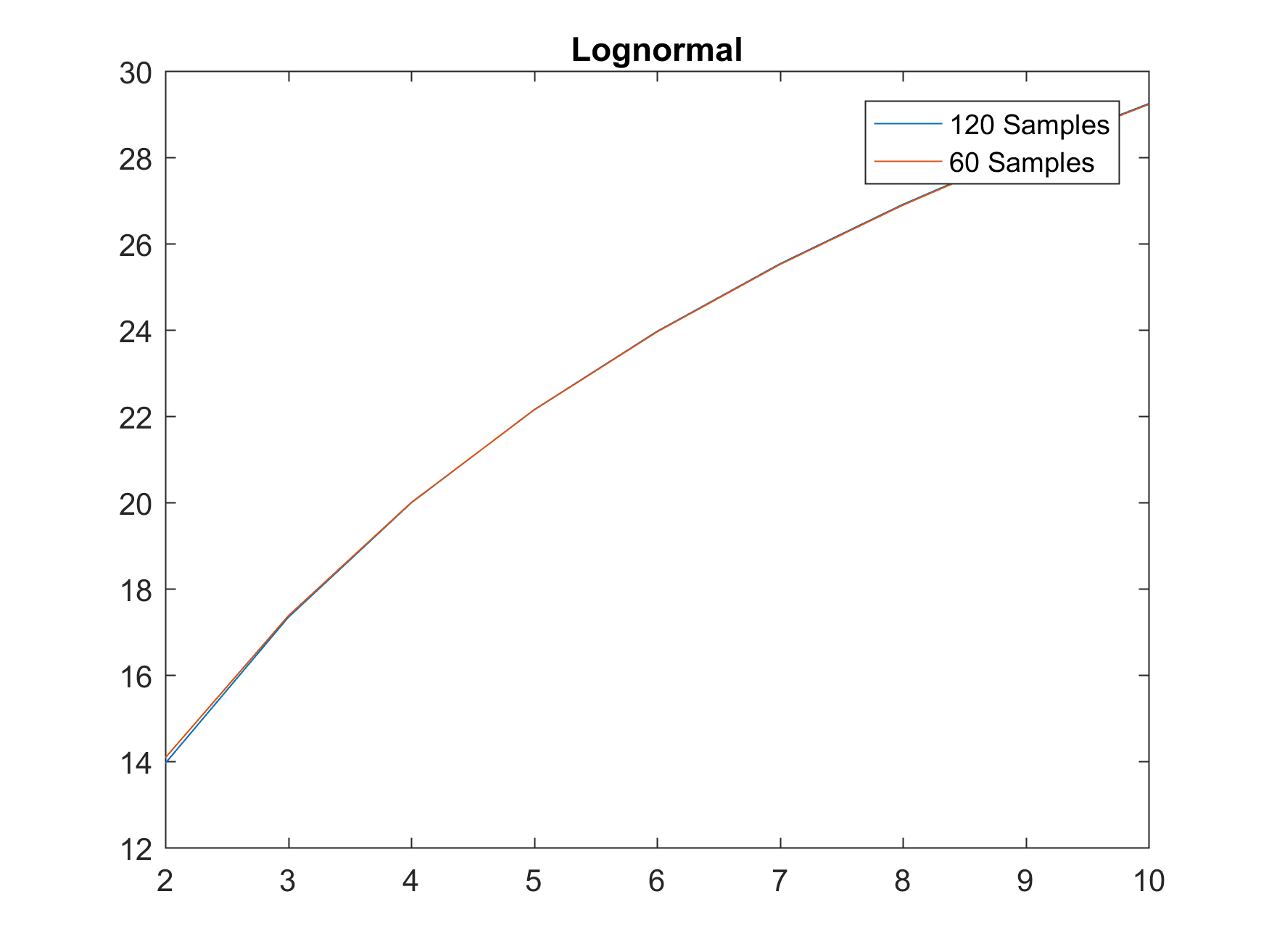}
	\caption{Divergence with different bandwidth parameter $B$. The randomness is log-normal distributed with dimension $m=11$.}
	\label{fig:bandwitdhlognrorm}
							\end{minipage}\vspace*{-2ex}%
			\end{figure}

 We investigate how the divergence between the reference and the true distributions varies with the bandwidth parameter used to estimate the reference. We consider a Gaussian distribution with dimension $m=11$, and sample sizes $N=120$ and $N=336$ (which are considered in Section \ref{sec:test_case_1}). Figures \ref{fig:bandwitdh120} and \ref{fig:bandwitdh336} show the KL divergence (estimated from 100,000 Monte Carlo samples drawn from the true distribution) against the bandwidth choice. In the figures we also show results with half of the samples sizes to give a sense of the sensitivity (and also motivated from the necessity of data splitting to be discussed momentarily). Among all the choices, $B=3$ appears the best as it gives the smallest divergence in three out of four different sample sizes. Figures \ref{fig:bandwitdh120d10} and  \ref{fig:bandwitdhlognrorm} further show the divergences between reference and true distributions when the truth follows other distributions, namely a Gaussian distribution with dimension $m=100$ and  a log-normal distribution with dimension $m=11$ respectively. We see that the graphs behave very differently from each other and the optimal bandwidth choices now deviate from 3, thus showing that the optimal bandwidth can depend heavily on the underlying distribution.

We note that the constructed $f_0$'s using kernel density estimation seem to be quite far from the true distribution. For example, in the problem considered in Section \ref{sec:test_case_1}, the KL divergence needs to be smaller than 1.25 in order to achieve a non-trivial solution. This is substantially smaller than 5.5, the lowest observed divergence value among all of Figures \ref{fig:bandwitdh120}--\ref{fig:bandwitdhlognrorm}. In other words, kernel density estimation is not efficient enough to obtain a good enough reference distribution for implementing DRO in this case. 

\subsubsection{Construction of Divergence-Based Uncertainty Sets}
Once we obtain a reference distribution, the next task is to calibrate the size of the uncertainty set. More precisely, we need to determine $\gamma$ for the set $\{f: \mathcal{D}(f\|f_0)\leq \gamma \}$, where $\mathcal D$ denotes the KL divergence, to cover the true distribution (with high confidence). This calls for the literature of divergence estimation. Here, we discuss the $k$-NN estimator studied by \cite{wang2009divergence,poczos2012nonparametric}. But before we proceed, we note that since $f_0$ itself is estimated from data, we need to be careful in controlling the statistical error in simultaneously estimating $f_0$ and $\gamma$. We consider two approaches. One is to use all the data to construct $f_0$ and reuse the same data to estimate $\gamma$. Another is to split the data into two groups, one for estimating $f_0$ and another for $\gamma$. In our experiments, the first approach turns out to consistently give a negative $\gamma$, indicating a poor estimation error (which is expected as the combined statistical error from $f_0$ and $\gamma$ is hard to control). Therefore, we adopt the second approach that splits the data.

\begin{figure}[t]
							\centering
							
							\begin{minipage}[b]{0.4\textwidth}
								\centering
								\includegraphics[width=\linewidth]{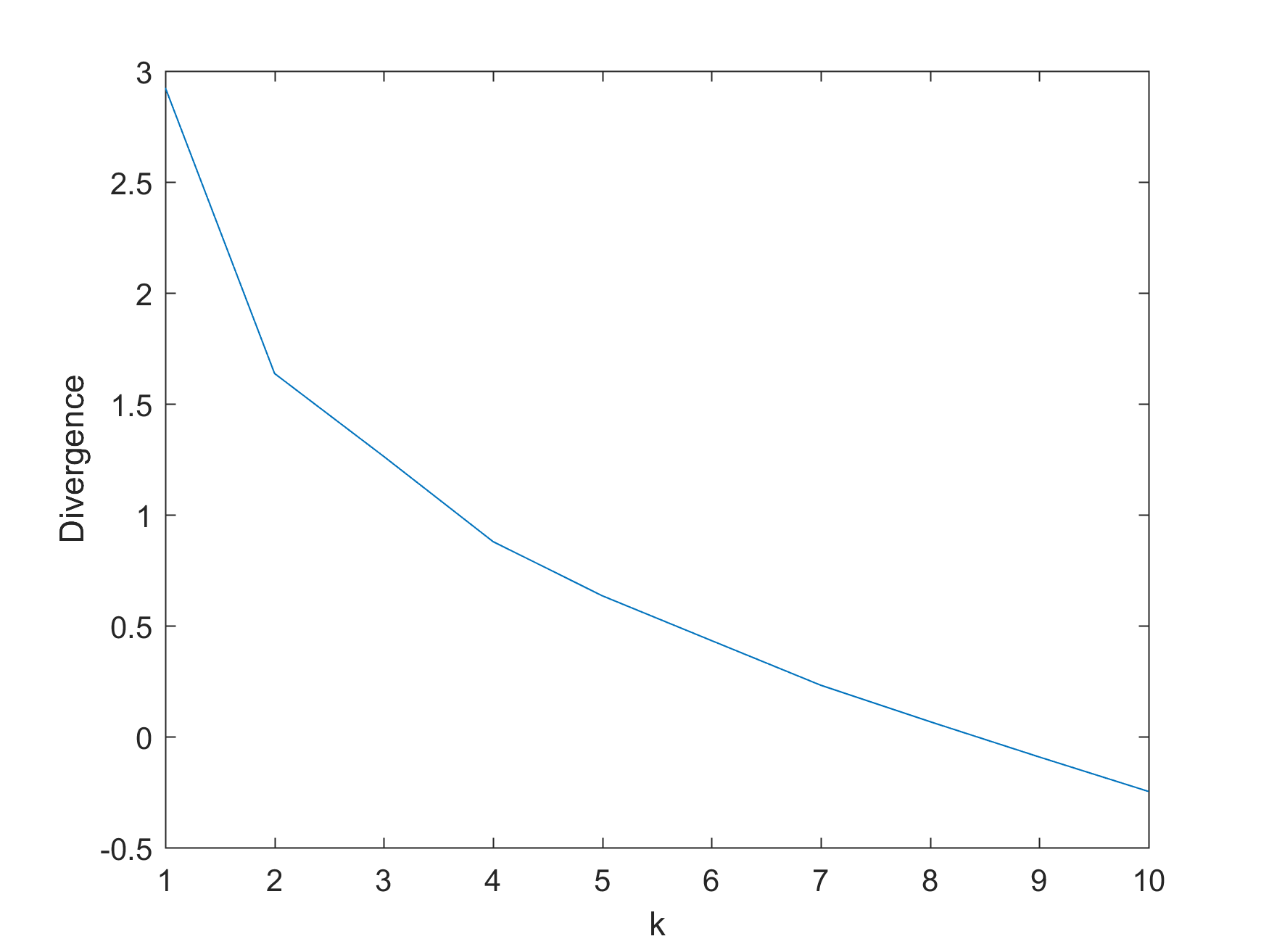}
								\caption{Estimated divergence with different $k$-NN parameter $k$ using 336 samples.}
								\label{fig:knn}
							\end{minipage}
                            \begin{minipage}[b]{0.4\textwidth}
								\centering
                            \includegraphics[width=\linewidth]{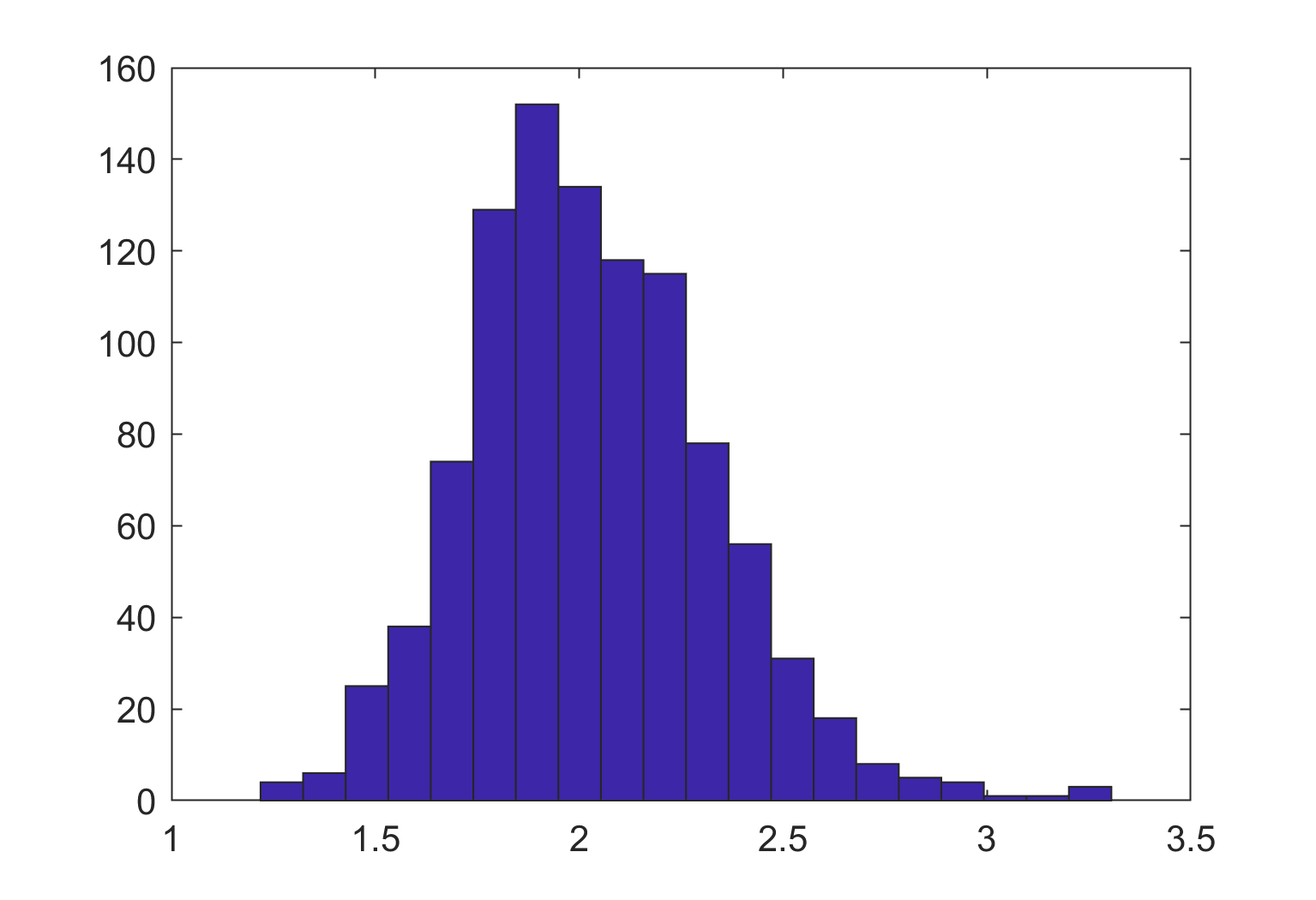}
	\caption{Histogram of the divergence estimates for $\gamma$ with $k=1$ and 336 samples.}
	\label{fig:bootmean}
							\end{minipage}\vspace*{-2ex}%
			\end{figure}

We investigate the quality of $k$-NN estimation with different choices of $k$, using an example of Gaussian distribution with $m=11$ and sample size 336. Here we split the data into two equal halves, and use the first half to estimate $f_0$ with bandwidth $B=3$ and the second half to estimate the divergence to calibrate $\gamma$. Figure  \ref{fig:knn} shows the average point estimate of the divergence using $k$-NN among 1000 experimental replications, against $k$. We see that  $k=1$ gives the closest estimate to the true divergence (5.5, using $B=3$ in Figure \ref{fig:bandwitdh336}). This observation is consistent with the known result in the literature that $k=1$ gives the smallest bias. However, even in this case the bias is still substantial, likely due to insufficient sample size. The performance is worse as $k$ increases.

Figure \ref{fig:bootmean} further shows the histogram of divergence estimates from 1000 experimental replications with $k=1$. The distribution of the estimates appears very spread out. Moreover, the biggest realized estimate (less than 3.5) is still far away from the true divergence (5.5 in Figure \ref{fig:bandwitdh336}). As noted in \cite{wang2009divergence}, estimating divergence for high-dimensional distributions with small sample typically incurs large variances and is challenging, in line with our observations here. For problems with even higher dimension (e.g., the setting in Figure \ref{fig:bandwitdh120d10}), we expect it to be even more difficult to obtain a reasonable divergence estimate.
			

\subsubsection{A Low-Dimensional Example for Divergence-Based DRO.}

\begin{table}[t]
	\centering
	\caption{Optimality and feasibility performances on a single $d=3$ dimensional linear CCP with Gaussian distribution for several methods, using different sample sizes. }

\begin{tabular}{|l||l|l|l|l|l|l||l|l|l|l||l|}
\hline
                     & RO                    & Recon           & SG          & FAST          & DRO Mo & DRO KL  &  SCA      \\ \hline
$n$                 & 120                    & 120                   & 120       & 120  & 120     & 10000          & -        \\ \hline
$n_1$                & 60                    & 60            & -            & 61       & -      & 5000        & -        \\ \hline
$n_2$               & 60                    & 60            & -            & 59                & -   & 5000        & -        \\ \hline
Obj. Val.           & -1000.62               & -1056.59       & -1062.54     & -1041.47   & -916.33  & -1036.19  &  -1051.08 \\ \hline
$\hat{\epsilon}$    & 0.0010  & 0.0166        & 0.0248      & 0.0149           & $3.52 \times 10^{-7}$ & 0.0055  & 0.0072        \\ \hline
$\hat{\delta}$      & 0                      & 0.023        & 0.004          & 0.003            & 0  & 0 & 0        \\ \hline
\end{tabular}

	\label{table.joint.normal.n=3}
\end{table}	

Since using divergence-based DRO in the previously considered problem size (e.g., dimension $11$) appears problematic, we investigate a very small problem with $d=3$. We consider a single linear CCP under a 3-dimensional Gaussian distribution. When we have only 120 samples, divergence-based DRO gives an adjusted tolerance level $\epsilon^*=8.5661e-07$, which is difficult to solve with sufficient accuracy using SAA and also likely leads to a very conservative solution. We therefore increase the sample size to 10,000, and finally obtain  $\epsilon^*= 0.0054$, which allows us to solve via SAA with 500 Monte Carlo samples. Table \ref{table.joint.normal.n=3} shows its performance and compares with other approaches. Divergence-based DRO (with sample size $10,000$) is now less conservative than moment-based DRO and our plain RO, but is outperformed by reconstructed RO, SG and FAST, all with a much smaller sample size (120), and SCA. 

\end{document}